\documentstyle[11pt,makeidx]{article}
\makeindex
\newcommand{\np}{\newpage}
\newcommand{\ROM}{\rm} %
\newcommand{\rmm}[1] {\rm{#1}}  
\newcommand{\bfm}[1] {\bf{#1}}           
\newcommand{\rmt}{\ROM\rm}
\newcommand{\itt}{\ROM\it}
\newcommand{\bft}{\ROM\bf}
\chardef\REST="16
\newcommand{\RESTR}{{\msas\REST}}
\newcommand{\restriction}{\RESTR}

\newcommand{\xieufm}{\font\teufm=eufm10 scaled \magstephalf
\font\seufm=eufm8 \font\zeufm=eufm8}
\newcommand{\xeufm}{\font\teufm=eufm10 
\font\seufm=eufm7 \font\zeufm=eufm7}
\newcommand{\ixeufm}{\font\teufm=eufm9 
\font\seufm=eufm7 \font\zeufm=eufm7}
\newcommand{\fraK}[1]{{\mathchoice {\hbox{\teufm{#1}}} 
{\hbox{\teufm{#1}}}{\hbox{\seufm{#1}}}{\hbox{\zeufm{#1}}} }}

\newcommand{\xieusm}{\font\teusm=eusm10 scaled \magstephalf
\font\seusm=eusm8 \font\zeusm=eusm8}
\newcommand{\xeusm}{\font\teusm=eusm10 
\font\seusm=eusm7 \font\zeusm=eusm7}
\newcommand{\ixeusm}{\font\teusm=eusm9 
\font\seusm=eusm7 \font\zeusm=eusm7}
\newcommand{\skr}[1]{{\mathchoice {\hbox{\teusm{#1}}} 
{\hbox{\teusm{#1}}}{\hbox{\seusm{#1}}}{\hbox{\zeusm{#1}}} }}

\newcommand{\ximsb}{\font\tmsbm=msbm10 scaled \magstephalf
\font\smsbm=msbm8 \font\zmsbm=msbm8}
\newcommand{\xmsb}{\font\tmsbm=msbm10 
\font\smsbm=msbm7 \font\zmsbm=msbm7}
\newcommand{\ixmsb}{\font\tmsbm=msbm9 
\font\smsbm=msbm7 \font\zmsbm=msbm7}
\newcommand{\BBB}[1]{{\mathchoice {\hbox{\tmsbm{#1}}} 
{\hbox{\tmsbm{#1}}}{\hbox{\smsbm{#1}}}{\hbox{\zmsbm{#1}}} }}

\newcommand{\ximsa}{\font\tmsam=msam10 scaled \magstephalf
\font\smsam=msam8 \font\zmsam=msam8}
\newcommand{\xmsa}{\font\tmsam=msam10 
\font\smsam=msam8 \font\zmsam=msam8}

\newcommand{\msas}[1]{{\mathchoice {\hbox{\tmsam{#1}}} 
{\hbox{\tmsam{#1}}}{\hbox{\smsam{#1}}}{\hbox{\zmsam{#1}}} }}
\def\addto#1#2{
\ifx\zone\undefined\let\zone=#1\def#1{\zone#2}\else
\ifx\ztwo\undefined\let\ztwo=#1\def#1{\ztwo#2}\else
\ifx\zthree\undefined\let\zthree=#1\def#1{\zthree#2}\else
\fi\fi\fi
}
\addto\normalsize   {\xieusm\ximsb\ximsa\xieufm}
\addto\small        {\xeusm\xmsb\xmsa\xeufm}
\addto\footnotesize {\ixmsb\ixeusm\ixeufm}
\newtheorem{gtheorem}{Theorem}
\newtheorem{theorem}{Theorem}
\newtheorem{assertion}[theorem]{Assertion}
\newtheorem{corollary}[theorem]{Corollary}
\newtheorem{definition}[theorem]{Definition}
\newtheorem{lemma}[theorem]{Lemma}
\newtheorem{proposition}[theorem]{Proposition}
\newtheorem{remark}[theorem]{Remark}
\newcommand{\TF}{\itt}
\newcommand{\bass}{\begin{assertion}\TF\ } 
\newcommand{\eass}{\end{assertion}}
\newcommand{\bcor}{\begin{corollary}\TF\ }
\newcommand{\ecor}{\end{corollary}}
\newcommand{\bdf} {\begin{definition}\rmt\ }
\newcommand{\edf} {\end{definition}} 
\newcommand{\ble} {\begin{lemma}\TF\ }
\newcommand{\ele} {\end{lemma}}
\newcommand{\bpro}{\begin{proposition}\TF\ } 
\newcommand{\epro}{\end{proposition}} 
\newcommand{\brem}{\begin{remark}\rmt\ }
\newcommand{\erem}{\end{remark}} 
\newcommand{\bte} {\begin{theorem}\TF\ }
\newcommand{\ete} {\end{theorem}}
\newcommand{\bgte} {\begin{gtheorem}\TF\ }
\newcommand{\egte} {\end{gtheorem}}

\newcommand{\proof}{\noi{\bft Proof}\hspace{3mm}}
\newcommand{\eproof}{\vspace{3mm}}
\newcommand{\ben}{\begin{enumerate}}
\newcommand{\een}{\end{enumerate}}
\newcommand{\bit}{\begin{itemize}}
\newcommand{\eit}{\end{itemize}}
\newcommand{\bce}{\begin{center}}
\newcommand{\ece}{\end{center}}
\newcommand{\bay}{\begin{array}}
\newcommand{\eay}{\end{array}}
\newcommand{\bqu}{\begin{quotation}}
\newcommand{\equ}{\end{quotation}}
\newcommand{\bbox}[1]{{\bfm{#1}}}
\newcommand{\BST} {\bbox{BST}}
  
\newcommand{\HST} {\bbox{HST}}
\newcommand{\IST} {\bbox{IST}}

\newcommand{\ZFC} {\bbox{ZFC}}
\newcommand{\aBI} {\bbox{BI}}

\newcommand{\aI}  {\bbox{I}}
\newcommand{\aS}  {\bbox{S}}
\newcommand{\aT}  {\bbox{T}}
\newcommand{\aB}  {\bbox{B}}

\newcommand{\ip}  {\bbox{IP}}
\newcommand{\ips}  {\ip^{\rbox{S}}}
\newcommand{\fa}{\bbox a}
\newcommand{\fb}{\bbox b}

\newcommand{\fx}{\bbox x}

\newcommand{\fA}{\bbox A}
\newcommand{\fB}{\bbox B}
\newcommand{\fG}{\bbox G}
\newcommand{\fL}{\bbox L}

\newcommand{\fT}{\bbox T}
\newcommand{\rbox}[1]{{\rmm{#1}}}
\newcommand{\tp}   {\tau}   
\newcommand{\card} {{\tt card}\,} 
\newcommand{\Ord}  {{\tt Ord}}
\newcommand{\Card} {{\tt Card}}
\newcommand{\sord} {{\tt SOrd}}
\newcommand{\dom}  {{\tt dom}\,} 
\newcommand{\ran}  {{\tt ran}\,} 
\newcommand{\Ind}  {{\tt Ind}} 
\newcommand{\st}   {{\tt st}\,}   
\newcommand{\stan}[1]  {{\tt Stan}\hspace{0.5pt}(#1)}
\renewcommand{\int}{{\tt int}\,} 
\newcommand{\Seq}  {{\tt Seq}}      
\newcommand{\Max}  {{\tt Max}\,}     
\newcommand{\Min}  {{\tt Min}\,}     
\newcommand{\irk}  {{\tt irk}\,} 
\newcommand{\nrk}  {{\tt nrk}\,} 
\newcommand{\sfo}  {\hspace{5pt}{\tt forc}\hspace{5pt}} 
\newcommand{\spfo}[1]{{\displaystyle
\hspace{5pt}{\tt forc}\hspace{1pt}_{#1}^{\phantom{C}}\hspace{5pt}}}
 
\newcommand{\rQ}{\rbox Q}
\newcommand{\al}{\alpha} 
\newcommand{\da}{\delta}

\newcommand{\ba}{\beta}
\newcommand{\ga}{\gamma}
\newcommand{\kpa}{\kappa}

\newcommand{\La}{\Lambda} 
\newcommand{\la}{\lambda} 
\newcommand{\sg}{\sigma}
\newcommand{\Sg}{\Sigma}
\newcommand{\vt}{\vartheta}
\newcommand{\vpi}{\varphi}

\newcommand{\za}{\zeta} 
\newcommand{\vep}{\varepsilon} 
\newcommand{\fsg}{{\cur\sg}}
\newcommand{\bbb}{\hspace{0.5pt}} 
\newcommand{\dvoj}[1]{\mathord{\bbb{\BBB #1}\bbb}}

\newcommand{\dE}{{\dvoj E}}
\newcommand{\dH}{{\dvoj H}}

\newcommand{\dI}{{\dvoj I}}

\newcommand{\dL}{{\dvoj L}}

\newcommand{\dP}{{\dvoj P}}
         \newcommand{\onep}{{\bf 1}\hspace{-0.5pt}_\dP}

\newcommand{\dN}{{\dvoj N}}
\newcommand{\dS}{{\dvoj S}}
\newcommand{\dV}{{\dvoj V}}
\newcommand{\cur}{}
\newcommand{\skrsp}{\hspace{0.5pt}} 
\newcommand{\scri}[1]{\mathord{\skrsp\skr #1\skrsp}} 
\newcommand{\cA}{{\scri A}}
\newcommand{\cC}{{\scri C}}
\newcommand{\cD}{{\scri D}}
\newcommand{\cE}{{\scri E}}

\newcommand{\cH}{{\scri H}}

\newcommand{\cL}{{\cal L}}
          \newcommand{\cLi}{{\cal L}^\infty}
          \newcommand{\Di}{D^\infty}
          \newcommand{\Dpi}{{D_+}^\infty}
          \newcommand{\Dsi}{\cD^\infty}

          \newcommand{\Epi}{{E_+}^\infty}
\newcommand{\cN}{{\scri N}}
\newcommand{\cP}{{\scri P}}
\newcommand{\cX}{{\scri X}}
\newcommand{\cY}{{\scri Y}}

\newcommand{\gotsp}{\hspace{0.5pt}}
\newcommand{\got}[1]{\mathord{\gotsp\fraK #1\gotsp}}  
\newcommand{\gA}{{\got A}}
\newcommand{\gB}{{\got B}}
\newcommand{\gC}{{\got C}}
\newcommand{\rh}{\rbox{h}}
\newcommand{\upr}{\mathord{\kern 1pt ^\rh\kern -1pt}}
\newcommand{\req}{\mathrel{\DS\kern1pt^\rh\kern-5pt=\kern1pt}}
\newcommand{\rin}{\mathrel{\DS\kern1pt^\rh\kern-5pt\in\kern 1pt}}
\newcommand{\rst}{{\DS{\upr\kern 0.5pt\st}}}
\renewcommand{\sup}[2]
{\mathord{\kern 0.05em\vphantom{X}^{#2}\kern -0.08em #1}}
\newcommand{\sups}[2]
{\mathord{\kern 0.05em\vphantom{X}^{#2}\kern -0.17em #1}}
\newcommand{\supi}[2]
{\mathord{\kern 0.05em\vphantom{X}^{#2}\kern -0.27em #1}}
\newcommand{\upsg}[1]{\sup{#1}\fsg}
\newcommand{\upsG}[1]{\sups{#1}\fsg}
\newcommand{\upSG}[1]{\supi{#1}\fsg}

\newcommand{\upa}[1]{\sups{#1}\ast}

\newcommand{\upst}{^{\rbox{st}}}
\newcommand{\upin}{^{\rbox{int}}}
\newcommand{\DS}{\displaystyle}
\newcommand{\Askip}{\hspace{0.5pt}}
\newcommand{\est}{{\DS\exists\upst}} 
\newcommand{\fst}{{\DS\forall\Askip\upst}}
\newcommand{\ein}{{\DS\exists\upin}}
\newcommand{\fin}{{\DS\forall\Askip\upin}}
\newcommand{\fstf}{{\DS\forall\Askip^{\rbox{stfin}}}}
\newcommand{\Sgs}{\Sg\upst}
\newcommand{\Ss}{{\DS{\Sgs_2}}}
\newcommand{\sq}{\subseteq}
\newcommand{\cj}{\mathbin{\hspace{2pt}\&\hspace{2pt}}}      
\newcommand{\orr}{\mathbin{\hspace{2pt}\textstyle\bigvee\hspace{2pt}}} 
\newcommand{\lra}{\longrightarrow} 
\newcommand{\llra}{\longleftrightarrow} 
\newcommand{\map}{\,\longmapsto\,} 
\newcommand{\res}{\mathbin{\restriction}}
\newcommand{\we}{{\mathbin{\kern 1.3pt ^\wedge}}}
\newcommand{\<}{\leq}

\newcommand{\ti}{\times}
\newcommand{\qed}{\hfill{$\msur\Box\msur$}}
\newcommand{\bult}{\hfill{$\msur\dashv\msur$}}
\newcommand{\dm}{$$}

\newcommand{\avl}[1]{|#1|}
\newcommand{\davl}[1]{\|#1\|}
\newcommand{\mo}{\models}
\newcommand{\ima} {\mathord{\hspace{0.1em}"\hspace{0pt}}}             
\newcommand{\emps}{\emptyset}
\newcommand{\fo}{\mathrel{{|\hspace{-1pt}|
\hspace{-3pt}\mathord{-}\hspace{-8pt}\mathord{-}} }}
\newcommand{\pfo}[1]{\mathrel{{\fo}_{#1}}}
\newcommand{\mem}{$\hspace{0.5pt}\msur\in\msur$}
\newcommand{\ste}[1]{\hbox{\hspace{0.5pt}{{\tt st}}-$\msur\in\msur$#1}}

\newcommand{\wfp}[2]{{\mathord{\ang{#1,#2}}}}
\newcommand{\ang} [1]{\langle #1\rangle}

\newcommand{\ans} [1]{\{\hspace{0.2mm}#1\hspace{0.2mm}\}}
\def\widehat{\overline}

\newcommand{\wx} {\widehat x}
\newcommand{\wy} {\widehat y}
\newcommand{\wS} {\widehat S}
\newcommand{\wla}{\widehat \la}

\newcommand{\brp}{{\breve p}}
       \newcommand{\brpi}{{\breve \pi}}
\newcommand{\brs}{{\breve s}}
\newcommand{\brx}{{\breve x}}
\newcommand{\bry}{{\breve y}}
\newcommand{\brG}{{\underline G}}

\newcommand{\eqg}{\mathbin{{\mathord=}_G}}
\newcommand{\ing}{\mathbin{{\mathord\in}_G}}
\newcommand{\ning}{\mathbin{{\mathord{\not\in}}_G}}
\newcommand{\inh}{\mathbin{{\mathord\in}_\dH}}
\newcommand{\ins}{\mathbin{{\mathord\in}_\dS}}
\newcommand{\ini}{\mathbin{{\mathord\in}_\dI}}
\newcommand{\inhp}{\mathbin{{\mathord\in}_{\dH'}}}
\newcommand{\stg}{{\st\!}_G\,}
\newcommand{\dd}[2]{\hbox{$\kern-0.7mm{#1}\kern-1mm$-#2}}
\newcommand{\noi}{\noindent}

\newcommand{\msur}{\protect\hspace{-1\mathsurround}}

\newcommand{\its}{\vspace{-1mm}}
\newcommand{\itla}{\item\label}
\newcommand{\ux} {x^0} 

\begin{document}
\normalsize
%
%

\title{Isomorphism property in nonstandard extensions of 
$\ZFC$ universe} 
\author{Vladimir Kanovei
\thanks
{\ 
\protect\begin{minipage}[t]{0.9\textwidth}
Moscow Transport Engineering Institute\protect\\
{\tt kanovei@mech.math.msu.su} \ and \ 
{\tt kanovei@math.uni-wuppertal.de}\protect\\
Partially supported by DFG grant and AMS grant 
\protect\end{minipage}\protect\vspace{2pt}
}
\and
Michael Reeken
\thanks
{\ Bergische Universit\"at -- 
GHS Wuppertal. \ {\tt reeken@math.uni-wuppertal.de}}
}
\date{01 March 1996}
\maketitle

\begin{abstract}
We study~\footnote
{\ The authors are pleased to mention useful conversations on 
the topic of the paper, personal and in written form, with
D.~Ballard, K.~Hrba\v cek, H.~J.~Keisler.  

The first (alphabetically) author is in debt to several 
institutions and personalities who facilitated his part of work 
during his visiting program in 1995 -- 1996, in particular 
University of Wuppertal, I.P.M. in Tehran, and 
Max Planck Institute at Bonn, and 
personally P.~Koepke and M.~J.~A.~Larijani.}     
models of $\HST$ 
(a nonstandard set theory which includes, in particular, the 
$\ZFC$ Replacement and Separation schemata 
and Saturation for well-order\-able families of internal sets). 

This theory admits an adequate formulation of the 
{\it isomorphism property\/} $\ip:$ any two elementarily 
equivalent internally presented structures of a well-orderable  
language are isomorphic. $\ip$ implies, for instance, that 
all infinite internal sets are equinumerous, and there exists 
a unique (up to isomorphism) internal elementary 
extension of the standard reals. 

We prove that $\ip$ is independent of $\HST$ (using the class of all 
sets constructible from internal sets) and consistent with $\HST$  
(using generic extensions of $\HST$ models by a sufficient number 
of generic isomorphisms).\\[3mm] 
\end{abstract}

\newpage

{\addtolength{\parskip}{-1pt}
\footnotesize\tableofcontents
\addtolength{\parskip}{1pt}}

\np

%
\section*{Introduction}
\addcontentsline{toc}{section}{\protect\numberline{}Introduction}

This article is a continuation of the authors' series of papers 
\cite{hyp1,hyp2,hyp3} devoted to set theoretic foundations of 
nonstandard mathematics. The aim of this paper is to 
accomodate~\footnote
{\ The question as how one can adequately develop advanced 
nonstandard tools like $\ip$ in the frameworks of the ``axiomatic'' 
setting of foundations of nonstandard analysis, that is, in a 
reasonable nonstandard set theory, was discussed in the course of 
a meeting between H.~J.~Keisler and one of the authors (V.~Kanovei, 
during his visit to Madison in December 1994).} 
an important nonstandard tool, the isomorphism 
property of Henson~\cite{he74}, in the context of an axiomatic 
treatment of nonstandard analysis. 

Let $\kpa$ be a cardinal in the $\ZFC$ universe. A nonstandard 
model is said to satisfy the 
\dd\kpa{}{\it isomorphism property\/}, $\ip_\kpa$ in brief, iff, 
\index{isomorphism property!ipk@$\ip_\kpa$}%
whenever $\cL$ is a first--order language of 
$\card \cL<\kpa,$ any two internally presented elementarily 
equivalent \dd\cL structures are isomorphic. 
(An \dd\cL structure $\gA=\ang{A;...}$ is {\it internally 
presented\/} if the base set $A$ and every interpretation under 
\index{internally presented structure}
$\gA$ of a symbol in $\cL$ are internal in the given nonstandard 
model.)  

Henson~\cite{he75}, Jin~\cite{j92,j92a}, Jin and Shelah~\cite{js} 
(see also Ross \cite{rs}) demonstrate that $\ip$ implies several 
strong consequences inavailable in the frameworks of ordinary 
postulates of nonstandard analysis, for instance the existence of 
a set of infinite Loeb outer measure which intersects every set of 
finite Loeb measure by a set of Loeb measure zero, the theorem that 
any two infinite internal sets have the same external cardinality, 
etc. 

In the course of this paper, we consider the following formulation 
of $\ip$ with respect to $\HST,$ 
\index{Hrbacek set theory@Hrba\v cek set theory, $\HST$}%
a nonstandard set theory in the 
\ste-language~\footnote
{\ The language containing $\in$ and $\st,$ the 
{\it standardness predicate\/}, as the atomic predicate symbols.}%
\index{language!stel@\ste-language}%
\index{standardness predicate@standardness predicate $\st$}%
, which reasonably models interactions between standard, internal, 
and external sets (see Section~\ref{s:hst}).

\bit
\item[$\ip:$] 
If $\cL$ is a first--order language containing 
(standard size)--many symbols then any two internally presented  
\index{isomorphism property!ip@$\ip$}
elementarily equivalent \dd\cL structures are isomorphic. 
\eit
(Formally, sets {\it of standard size\/} 
\index{set!of standard size}
are those equinumerous to a set of the 
form: $\upsG S=\ans{x\in S:\st x},$ where $S$ is standard and 
\index{zzsigma@$\upsG X$}%
\index{zzstx@$\st x,$ the standardness predicate}%
\index{set!standard set}%
$\st x$ means: $x$ is a standard set.~\footnote
{\ Take notice that in $\HST$ {\it standard size\/} is equivalent to 
each of the following: {\it wellorderable\/}, 
{\it equinumerous to a wellfounded set\/}, 
{\it equinumerous to an ordinal\/}.}%
)

The following is the main result of the paper referred to in the title.

\bte
\label{maint}
$\ip$ is consistent with and independent of\/ $\HST.$ 
In addition, let\/ $\fT$ be a theory\/ $\HST + \ip$ or\/ 
$\HST + \neg\;\ip.$ Then\its
\ben 
\def\theenumi{(\Roman{enumi})}
\def\labelenumi{{\rmt\theenumi}}
\itla{esis}\msur
$\fT$ is equiconsistent with $\ZFC$.\its

\itla{serv}\msur
$\fT$ is a conservative extension of\/ $\ZFC$ in the following 
sense. Let\/ $\Phi$ be a closed\/ \mem-formula, $\Phi\upst$ be 
the formal relativization to the predicate\/ $\st.$ Then\/ 
\index{zzphist@$\Phi\upst$}%
$\Phi$ is a theorem of\/ $\ZFC$ iff\/ $\Phi\upst$ is a theorem 
of\/ $\fT$.~\footnote
{\rmt\ In other words it is asserted that $\ZFC$ proves $\Phi$  
iff $\HST$ proves that $\Phi$ holds in the standard subuniverse.}
\its

\itla{me}
Every countable model\/ $\dS\mo\ZFC$ can be embedded, as the 
class of all standard sets, in a model\/ $\dH$ of $\fT,$ 
satisfying the following additional property \ {\rm\ref{red}$:$}\its

\itla{red}
If\/ $\Phi(x_1,...,x_n)$ is a\/ \ste-formula then there exists 
an\/ \mem-formula\/ $\Phi^\ast(x_1,...,x_n)$ such that, for 
all sets\/ $x_1,...,x_n\in\dS,$ $\Phi(x_1,...,x_n)$ is true in\/ 
$\dH$ iff $\Phi^\ast(x_1,...,x_n)$ is true in\/ $\dS$.~\footnote
{\rmt\ Thus the truth of \ste-formulas with standard 
parameters in $\dH$ can be investigated in $\dS.$ A similar 
property was called {\it Reduction\/} in \cite{hyp1}.}
\een 
\ete
%
%
%
The $\HST$ models involved in the proof of Theorem~\ref{maint} are 
obtained as the results of several consecutive extensions of an 
initial model $\dS$ of $\ZFC;$ $\dS$ becomes the class of all 
\index{class!s@$\dS$ of all standard sets}
standard sets in the final and intermediate models. The sequence 
of extensions contains the following steps: 

{\it Step 1\/}. 
We extend $\dS$ to a model $\dS^+$of $\ZFC$ plus global 
choice, adjoining a generic wellordering of the universe by an 
old known method of Felgner~\cite{fe}. $\dS$ and $\dS^+$ 
contain the same sets. 

{\it Step 2\/}. 
We extend $\dS$ to a model $\dI$ of {\it bounded set thery\/} 
\index{class!i@$\dI$ of all internal sets}
$\BST,$ a nonstandard set theory similar to $\IST$ of 
\index{bounded set theory@bounded set theory, $\BST$}%
Nelson~\footnote
{\ It is not known whether $\IST$ itself admits the treatment 
similar to steps 1 -- 6.}
, using a global choice function from $\dS^+$ to define $\dI$ 
as a kind of ultrapower of $\dS.$ $\dS$ is the class of all 
standard sets in $\dI.$ This step was described in \cite{hyp1}.

{\it Step 3\/}. We extend $\dI$ to a model $\dH$ of 
\index{class!h@$\dH$ of all external sets}
{\it Hrba\v cek set theory\/} $\HST,$ a nonstandard set theory 
\index{Hrbacek set theory@Hrba\v cek set theory, $\HST$}%
which contains, for instance, Separation and Collection in the 
\ste-language, and Saturation for standard size families of 
internal sets. The universe $\dH$ is isomorphic to an inner 
\ste-definable structure in~$\dI$. 
This step was described in~\cite{hyp2}. Elements 
of $\dH$ are essentially those sets which can be obtained from 
sets in $\dI$ by the procedure of assembling sets along 
wellfounded trees definable in $\dI.$ $\dS$ 
remains the class of all standard sets in $\dH$ while $\dI$ 
becomes the class of all {\it elements\/} of standard sets (that 
is, {\it internal\/} sets by the formal definition) in $\dH$.  
 
We proved in \cite{hyp1,hyp2} that $\BST$ and $\HST$ 
are equiconsistent with $\ZFC$ and are conservative extensions 
of $\ZFC$ in the sense of statement~\ref{serv}. 

{\it Step 4\/}. Given a model $\dH$ of $\HST,$ we define, in  
Section~\ref{str}, the class $\dL[\dI]$ 
\index{class!li@$\dL[\dI]$ of all sets constructible from internal 
sets}%
of all sets constructible in 
$\dH$ from internal sets (an inner class in $\dH$). The particular 
case we consider (constructibility from internal sets) allows 
to define the constructible sets much easier than in $\ZFC,$ 
because $\dI$ contains all ordinals and essentially all  
patterns of constructions which might be involved in definitions  
of constructible sets; this allows to avoid any kind of transfinite 
recursion in the process.

We prove that $\dL[\dI]$ is a model of $\HST$ which satisfies 
certain additional properties, in particular \dd\dI infinite 
internal sets of different \dd\dI cardinalities remain 
non--equinumerous in $\dL[\dI],$ so that $\dL[\dI]$ models the 
negation of $\ip.$ 
This leads to the proof of different parts of the 
theorem, with respect to the theory $\HST+\neg\;\ip$. 

We prove that in addition $\dL[\dI]$ satisfies the following 
choice--like statement: for any cardinal $\kpa,$ every 
\dd\kpa closed p.\ o.\ set is 
\index{set!kclosed@\dd\kpa closed}%
\index{set!kdistributive@\dd\kpa distributive}%
\dd\kpa distributive, which is of great importance for the 
further use of $\dL[\dI]$ as the ground model for generic  
extensions.

{\it Step 5\/} (actually irrelevant to the content of this paper). 
Given a standard cardinal $\kpa,$ $\HST$ admits a subuniverse 
$\dH_\kpa$ (an inner \ste-definable class in $\dH$) which models a 
\dd\kpa version of $\HST$ (the Saturation and standard size Choice 
suitably restricted by $\kpa$) plus the Power Set axiom. The class 
$\dI_\kpa=\dI\cap\dH_\kpa$ of all internal elements in $\dH_\kpa$ 
is equal to the collection of all $x\in\dI$ which belong to a 
standard set of \dd\dS cardinality $\<\kpa.$ There also exist 
subuniverses $\dH_\kpa'$ which model a weaker \dd\kpa version of 
$\HST$ plus Power Set and full Choice. These subuniverses 
are defined and studied in \cite{hyp3}.

{\it Step 6\/}. To get a model for $\HST+\ip,$ we construct a 
generic extension $\dL[\dI][G]$ of an $\HST$ model of the form 
$\dL[\dI]$ (or any other model $\dH$ of $\HST$ satisfying the 
abovementioned choice--like principle), where $G$ is a generic 
class, essentially a collection of generic isomorphisms between 
suitable internally presented elementarily equivalent structures. 

We show how to develop forcing in $\HST$ in Section~\ref{f}. The 
forcing technique in principle resembles the classical $\ZFC$ 
patterns. However there are important differences. In particular 
to prevent appearance of new collections of standard sets in 
generic extensions (which would contradict {\it Standardization\/}, 
one of the axioms of $\HST$), the forcing notion must be 
{\it standard size distributive\/}. More differences appear in the 
{\it class\/} version of forcing introduced in Section~\ref{total}. 
In particular, to provide Separation and Collection in the class 
generic extensions we consider, a permutation technique is used; in 
the $\ZFC$ version, this tool is usually applied for different 
purposes. 

We demonstrate in Section~\ref{isom} how to force a generic 
isomorphism between two particular internally presented 
elementarily equivalent structures $\gA=\ang{A;...}$ and 
$\gB=\ang{B;...}$ in a model $\dH$ of $\HST.$ The forcing notion 
consists of internal $1-1$ maps from an internal $A'\sq A$ onto 
an internal $B'\sq B,$ such that every $a\in A'$ behaves in $\gA$ 
completely as $b=p(a)$ behaves in $\gB.$ This requirement means, 
for instance, that $a$ 
satisfies an \dd\cL formula $\Phi(a)$ in $\gA$ iff $b=p(a)$ 
satisfies $\Phi(b)$ in $\gB.$ But not only this. The main 
technical problem is how to expand a condition $p$ on some 
$a\in A\setminus \dom p,$ in other words, to find an appropriate 
counterpart $b\in B\setminus \ran p$ which can be taken as $p(a).$ 
To carry out this operation, we require that conditions $p$ 
preserve sentences of a certain type--theoretic extension of $\cL,$ 
the original first--order language, rather than merely 
\dd\cL formulas.~\footnote
{\ This type of forcing may lead to results related to 
nonstandard models in the $\ZFC$ universe, like the following: 
if $U$ is an \dd{\aleph_1}saturated nonstandard structure then 
there exists a generic extension of the $\ZFC$ universe where  
$U$ remains \dd{\aleph_1}saturated and satisfies 
$\ip_{\aleph_1},$ in the sense of locally internal isomorphisms. 
But this is a different story.}

It is worth noticing that the generic isomorphisms $H$ obtained by 
this forcing satisfy an interesting additional requirement. They 
are {\it locally internal\/} in the sense that, unless $A$ and $B$ 
are sets of standard finite number of elements, for any $a\in A$ 
there exists an internal set $A'\sq A,$ containing more than a 
standard finite number of elements, such that $H\res A'$ is 
internal. 

Section~\ref{total} demonstrates how to gather different generic 
isomorphisms in a single generic class using product forcing 
with internal \dd\dI finite support. (Fortunately new internally 
presented structures do not appear, so that the product rather than 
iterated forcing can be used here.) This results in a theorem which 
says that every countable model of $\HST$ satisfying the 
abovementioned additional property admits an extension with the 
same standard and internal sets, which satisfies $\HST+\ip.$ This 
leads to the proof of different parts of Theorem~\ref{maint}, with 
respect to the theory $\HST+\ip$. 

The countability assumption is used here simply as a sufficient 
condition for the existence of generic sets. (We shall in fact, 
for the sake of convenience, consider {\it wellfounded\/} $\HST$ 
models --- those having a well-founded class of ordinals in the 
wider universe --- but show how one obtains the result in the 
general case.) If the ground model is not assumed to be countable, 
a  Boolean--valued 
extension is possible, but we shall not proceed in this way. 

Section~\ref{final} completes the proof of Theorem~\ref{maint}. 

We use the ordinary set theoretic notation, with perhaps one 
exception: the \dd f{\it image\/} of a set $X$ will be denoted by 
\index{image@image $f\ima X$}%
\index{zzfimax@$f\ima X$}%
$f\ima X=\ans{f(x):x\in X}.$ The model theoretic notation will be 
elementary, self-explanatory, and consistent with Chang and 
Keisler~\cite{ck}. The reader is assumed to have an 
acquaintance with forcing and basic ideas and technique of 
nonstandard mathematics.

\subsubsection*{Remark}


It is sometimes put as 
a reservation by opponents of the nonstandard approach that the 
nonstandard real numbers are not uniquely defined, as long as one 
defines them by different constructions like ultrapowers of the 
``standard'' reals. (See Keisler~\cite{keis??} for a discussion of 
this matter.) 

It is obvious that any typical nonstandard set theory defines the 
standard reals and a nonstandard extension of the real line (the 
reals in the internal universe) uniquely. However $\HST$ does a 
little bit more: it is a particular property of this theory that 
the internal universe $\dI$ admits an \mem-definition in the 
external (bigger) universe $\dH,$ as a class of all sets $x$ 
such that there exists $y$ satisfying the property that the set 
$\ans{z:x\in z\in y}$ is linearly ordered but not well-founded. 
Thus the nonstandard reals are in a sense \mem-{\it unique\/}, not 
merely \ste-{\it unique\/}, in $\HST.$ (Without important changes 
in the set-up a result like this is hardly 
possible in nonstandard ``superstructures'', where one drops all 
the memberships from the ground level to carry out the 
construction.) 

The isomorphism property $\ip,$ if it holds in the external 
universe, makes the uniqueness much more strong: simply all 
internally presented elementary 
extensions of the standard reals are mutually isomorphic. 

As for the standard reals, they are also \mem-unique in $\HST,$ 
up to isomorphism at least, because the standard subuniverse 
$\dS$ is isomorphic in $\HST$ to an \mem-definable class, the 
class $\dV$ of all wellfounded sets. (This 
propertry in principle allows to develop mathematics in $\HST$ 
in terms of asterisks rather than the standardness predicate, 
see Subsection~\ref{change} below.)

On the other hand, Theorem~\ref{maint} makes it clear that, as 
long as one is interested in the study of standard mathematical 
objects, one can legitimately consider things so that the standard 
universe $\dS$ is in fact the standard part of a wider universe 
$\dH$ of $\HST+\ip,$ where different phenomena of ``nonstandard'' 
mathematics can be adequately presented.

\np

%
\section{Basic set theory in $\HST$}
\label{s:hst}

The development of $\HST$ in this section is based in part on 
ideas in early papers on external set theory of Hrba\v cek 
\cite{h78,h79} and Kawa\"\i~\cite{kaw83}.

\subsection{The axioms}
\label{hst:ax}

Hrba\v cek set theory $\HST$ (introduced in \cite{hyp2} on the base 
\index{Hrbacek set theory@Hrba\v cek set theory, $\HST$}
of an earlier version of Hrba\v cek \cite{h78}) is a theory in the \ste-language. 
It deals with three types of sets, standard, internal, and external. 
\index{set!standard set}
\index{set!internal set}
\index{set!external set} 

{\it Standard\/} sets are those $x$ satisfying $\st x.$ 
{\it Internal\/} sets are those sets $x$ which satisfy $\int x,$ 
where $\int x$ is the \ste-formula $\est y\;(x\in y)$ 
\index{zzintx@$\int x$}
(saying: $x$ belongs to a standard set). Thus the internal sets are 
precisely all the elements of standard sets. 
{\it External\/} sets are simply all sets.

\bdf
\label{his}
$\dS,\;\,\dI,\;\,\dH$ will denote the classes of all 
\index{class!s@$\dS$ of all standard sets}%
\index{class!i@$\dI$ of all internal sets}%
\index{class!h@$\dH$ of all external sets}%
standard and all internal sets, and the universe of all sets, 
respectively. 

\index{zzsigma@$\upsG X$}
$\upsG X=\ans{x\in X:\st x}=X\cap\dS$ for any set $X$.\qed
\edf

\newcommand{\iits}{\vspace{0pt}}

$\HST$ includes the following axioms:\iits
\ben
\def\theenumi{Ax.\arabic{enumi}}
\def\labelenumi{\theenumi}
\itla{it1} 
Axioms for standard and internal sets:
\ben 
\itla{zfc} 
$\msur\Phi\upst,$ where $\Phi$ is an arbitrary axiom of 
$\ZFC$
;\vspace{1mm}

\itla{tfer} 
{\it Transfer}:
\index{axiom!Transfer}
$\;\ein\,x\;\Phi\upin (x)\;\;\lra\;\;\est x\;\Phi\upin(x)$,\\[1mm]
where $\Phi$ is an \mem-formula containing only standard 
parameters, and $\Phi\upin$ denotes relativization of $\Phi$ 
to the class $\ans{x:\int x}$;\vspace{1mm}

\itla{trai} 
$\msur\fin x\;\forall\,y\in x\;(\int y)$: 
\ transitivity of the internal subuniverse.
%
\een

\itla{sza} 
{\it Standardization}:  
\index{axiom!Standardization}
$\forall\,X\;\est Y\;(\upsG X=\upsg Y)$.\iits 

\itla{hzfc} 
The $\ZFC$ Pair, Union, Extensionality, Infinity axioms, 
together with Separation, Collection, Replacement for all 
\ste-formulas.\iits

\itla{wr} 
{\it Weak Regularity\/}: 
\index{axiom!Weak Regularity}
for any nonempty set $X$ there exists 
$x\in X$ such that $x\cap X$ contains only internal elements.\iits 

\itla{sat} 
{\it Saturation\/}: 
\index{axiom!Saturation}
if $\cX$ is a set of standard size such that every $X\in \cX$ is 
internal and the intersection\/ $\bigcap \cX'$ is nonempty for any 
finite nonempty $\cX'\sq \cX,$ then $\bigcap \cX$ is nonempty.\iits

\itla{cho} 
\index{axiom!standard size Choice}
Choice in the case when the domain $X$ of the choice 
function is a set of standard size ({\it standard size Choice\/}), 
and Dependent Choice.\iits
\een

\subsection{Comments on the axioms}
\label{hst:comm}

The quantifiers $\est,\;\fst,\;\ein$ above have the obvious 
\index{zzqest@$\est$}
\index{zzqfst@$\fst$}
\index{zzqein@$\ein$}
meaning (there exists a standard ... , etc.). The first two will 
be of frequent use below. 

Axiom schema \ref{zfc} says that $\dS,$ the class of all 
standard sets, models $\ZFC.$ (Of course the $\ZFC$ 
Separation, Collection, and Replacement schemata are assumed to 
be formulated in the \mem-language in this item.) As an 
immediate corollary, we note that this implies $\dS\sq\dI.$ 

Transfer \ref{tfer} postulates that $\dI,$ 
the universe of all internal sets, is an 
elementary extension of $\dS$ in the \mem-language. 
Axiom \ref{trai} says that the internal sets form the ground 
in the \mem-hierarchy in $\dH,$ the main universe. 

We do not include here $\BST\upin,$ all the axioms of $\BST$ 
(see Subsection~\ref{bst:in:i} on the matters of $\BST$) 
relativized to the subuniverse $\dI$ of all internal sets, to the 
list of axioms, as it was the case in \cite{hyp2,hyp3}. Only 
Transfer and $\ZFC$ are included explicitly. But the rest of the 
$\BST$ axioms are more or less simple corollaries of other axioms, 
see Proposition~\ref{bstbi}. 

Standardization \ref{sza} is very important: it guarantees 
that $\dH$ does not contain collections of standard sets other 
than those which essentially already exist in $\dS.$ A simple 
corollary: a set in $\dH$ cannot contain all standard sets. One 
more application is worth to be mentioned.

\ble
\label{hstb}
{\rm [\hspace{1pt}Boundedness\hspace{1pt}]} \  
If\/ $X\sq\dI$ then\/ $X\sq S$ for a standard $S$.
\ele
\proof Each $x\in X$ is internal, hence belongs to a standard $s.$ 
By the Collection axiom, there is a set $B$ such that every 
$x\in X$ belongs to a standard $s\in B.$ By Standardization, there 
exists a standard set $A$ containing the same 
standard elements as $B$ does. We put $S=\bigcup A$.\qed\eproof

Group \ref{hzfc} misses the Power Set, Choice, and 
Regularity axioms of $\ZFC.$ Choice and Regularity still 
are added in weaker forms below. This is not a sort 
of incompleteness of the system; in fact each of the three 
mentioned axioms contradicts $\HST$. 

Axiom \ref{wr} says that all sets are well-founded over the 
universe $\dI$ of internal sets, in the same way as in $\ZFC$ all 
sets are well-founded over $\emps.$ (Take notice that $\dI$ itself 
is \underline{not} well-founded in $\dH:$ e.\ g.\ the set of all 
nonstandard 
\dd\dI natural numbers does not contain an \mem-minimal element.) 
There is, indeed, an essential difference with the $\ZFC$ 
setting: now $\dI,$ the ground level, explicitly contains a 
sufficient amount of information about the ordinals which 
determine the cumulative construction of $\dH$ from $\dI$.

{\it Sets of standard size\/} are those of the form 
\index{set!set of standard size}
$\ans{f(x):x\in \upsG X},$ where $X$ is standard and $f$ 
any function. However we shall see that in $\HST,$ 
``standard size'' = ``well-orderable'' = 
``equinumerous to a well-founded set''. 
The notion of a finite set in Axiom~\ref{sat} will be commented 
upon below.

It was convenient in \cite{hyp2,hyp3} to include one more 
axiom, {\it Extension\/}, to the list of $\HST$ 
axioms. Here we obtain it as a corollary. 

\ble
\label{exten}
{\rm [\hspace{1pt}Extension\hspace{1pt}]} \  
Suppose that\/ $S$ is a standard set and\/ $F$ is a function 
defined on the set\/ ${\upsG S=\ans{x\in S:\st x}},$ and\/ $F(x)$ 
contains internal elements for all\/ ${x\in \upsG S}.$ Then there 
exists an {\it internal\/} function\/ $f$ defined on\/ $S$ and 
satisfying\/ $f(x)\in F(x)$ for every\/ $x\in\upsG S$.
\ele
\proof We use the standard size Choice to obtain a (perhaps 
non--internal) function $g:\upsG S\,\lra\,\dI$ satisfying 
$g(x)\in F(x)$ for all standard $x\in S.$ It remains to apply 
Saturation to the family of (obviously internal) sets 
$G_x=\ans{f\in\dI:\dom f=S\cj g(x)=f(x)},$ where $x\in\upsG S$.
\qed

\subsection{Condensation of standard sets to well-founded sets}
\label{conden}

It is a typical property of nonstandard structures that 
standard sets reflect a part of the external universe. In the $\HST$ 
setting, this phenomenon appears in an isomorphism between the 
standard subuniverse $\dS$ and a certain transitive subclass of 
$\dH$ -- the class of all well-founded sets. 

\bdf
\label{dfcond}
We define $\wx=\ans{\wy:y\in\upsG x}$ for every $x\in\dS$.
\index{zzxhat@$\wx$}
 
We set $\dV=\ans{\wx:x\in\dS}$ (the {\it condensed 
subuniverse\/}.)\qed
\index{class!v@$\dV$ of all well-founded sets}
\edf
The next lemma shows that this definition is legitimate in $\HST$. 

\ble
\label{wfS}
The restriction\/ ${\in}\res\dS$ is a well-founded relation in $\dH$.
\ele
\proof Consider a nonempty set $X\sq\dS.$ By Standardization, there 
exists a standard set $S$ such that $X=\upsG S=S\cap\dS.$ Since 
$\dS$ models $\ZFC,$ $S$ contains, in $\dS,$ an \mem-mimimal 
element $s\in S.$ Then $s$ is \mem-minimal in $S$ also in the 
subuniverse $\dI$ 
by Transfer, and in $\dH$ by the definition of $\dI$.\qed\eproof

Thus $\wx$ is well defined in $\dH$ for all $x\in\dS$.

\ble
\label{s2w}
The map\/ $x\map\wx$ is an\/ \mem-isomophism\/ $\dS$ onto\/ $\dV.$ 
$\dV$ is a transitive class in\/ $\dH,$ and an\/ \mem-model of\/ 
$\ZFC.$ 
Every subset of\/ $\dV$ belongs to $\dV$.
\ele
\proof We have to prove first that $x=y$ iff $\wx=\wy,$ and $x\in y$ 
iff $\wx\in\wy,$ for all $x,\,y\in\dS.$ Only the direction 
$\longleftarrow$ is not obvious. Let $\wx=\wy.$ Then for 
each standard $x'\in x$ there exists standard $y'\in y$ such that 
$\wx'=\wy',$ and vice versa. This observation, plus Transfer (to 
see that standard sets having the same standard elements are equal) 
provides the proof of the first assertion by the induction on the 
\mem-rank of $x,\,y$ in $\dS$ (based on Lemma~\ref{wfS}). 

To prove the second assertion, let $\wx\in\wy.$ Then by definition 
$\wx=\wx'$ for some standard $x'\in y.$ Therefore $x=x'\in y,$ 
as required.

Let $W\sq \dV.$ Then $X=\ans{x:\wx\in W}$ is a subset of $\dS,$ so 
$X=\upsG S$ for a standard $S$ by Standardization. It follows that 
$W=\wS,$ as required.\qed\eproof 
 
Thus we have a convenient transitive copy $\dV$ of $\dS$ in $\dH$. 

Let a {\it well-founded set\/} mean: a set $x$ which belongs to 
\index{set!well-founded set}
a transitive set $X$ such that ${\in}\res X$ is a 
well-founded relation. (Axioms of group~\ref{hzfc} suffice to 
prove that every set belongs to a transitive set $X,$ but the 
membership on $X$ may be ill-founded. Take e. g.  
the set of all \dd\dI natural numbers.) 

\ble
\label{Wwf}
$\dV$ is the class of all well-founded sets in $\dH$.
\ele 
\proof Every $w\in\dV$ is well-founded in $\dH$ by Lemma~\ref{s2w}, 
because this is true for sets in $\dS.$ Suppose that $W\in\dH$ is 
well-founded, and prove that $W$ belongs to $\dV.$ We can assume 
that $W$ is transitive and ${\in}\res W$ is a well-founded relation. 
In this assumption, let us prove that $w\in\dV$ for all $w\in W,$ by 
\mem-induction. Suppose that $w\in W,$ and it is already known 
that each $w'\in w$ belongs to $\dV.$ Then $w\in\dV$ by 
Lemma~\ref{s2w}.\qed

\bpro
\label{memdef}
In\/ $\dH,$ $x\in\dI$ iff there is a set\/ $y$ such that the 
``interval''\/ $\ans{z:x\in z\in y}$ is linearly ordered by\/ 
$\in$ but not well-ordered.\qed
\epro
(This will not be used below, so we leave the proof for the 
reader.)

\subsection{Ordinals and cardinals in the external universe}
\label{orcar}

$\ZFC$ admits several formally different but equivalent definitions 
of ordinals. Since not all of them remain equivalent in $\HST,$ let 
us make it clear that an {\it ordinal\/} is a transitive set 
\index{ordinal} 
well-ordered by $\in.$ The following lemma will be used to prove 
that $\dH$ and $\dV$ contain the same ordinals.

\ble
\label{wo=ss}
Every set\/ $w\in \dV$ can be well-ordered and has standard size in\/ 
$\dH.$ Conversely if\/ $z\in\dH$ is a set of standard size or can 
be well-ordered in\/ $\dH$ then\/ $z$ is equinumerous with some 
$w\in \dV$ in $\dH$.
\ele
\proof Each set $w\in\dV$ is a set of standard size in $\dH$ 
because we have $w=\wx=\ans{\wy:y\in\upsG x}$ for a standard 
$x.$ To prove that $w$ can be well-ordered, it suffices to check 
that $\upsG x=x\cap\dS$ can be well-ordered in $\dH.$ Let $<$ be a 
standard well-ordering of $x$ in $\dS.$ Then $<$ may be not a 
well-ordering of $x$ in $\dH$ since $x$ obtains new 
subsets in $\dH.$ But $<$ still well-orders $\upsG x.$ (Indeed 
let $u'\sq \upsG x$ be a nonempty set. Then $u'=\upsG u$ for a 
standard set $u\sq x$ by Standardization. Take the \dd<least 
element of $u$ in $\dS$.)

To prove the converse, note that, in $\dH,$ every set of 
standard size can be well-ordered --- by the previous argument. 
Thus let $Z\in\dH$ be well-ordered in $\dH;$ let us 
check that $Z$ is equinumerous with a set $W\in\dV$.

Since the class $\sord$ of all \dd\dS ordinals (i. e. standard 
sets which are ordinals in $\dS$) is well-ordered by Lemma~\ref{wfS}, 
either there exists an order preserving map: $\sord$ onto an initial 
segment of $Z$ or there exists an order preserving map: $Z$ onto a  
proper initial segment of $\sord$. 

The ``either'' case is impossible by axioms~\ref{hzfc} and 
Lemma~\ref{hstb}. Thus we have the ``or'' case. Let $\la$ be the 
least standard ordinal which does not belong to the proper initial 
segment. Then we have a $1-1$ map from $\upsG \la$ onto $Z.$ Then 
$W=\wla\in\dV$ admits a $1-1$ map onto $Z,$ as required.\qed

\bcor
\label{o}
Universes\/ $\dH$ and\/ $\dV$ contain the same ordinals.
\ecor
\proof 
Suppose that $\xi\in\dV$ is an ordinal in $\dV.$ Then $\xi$ 
remains an ordinal in $\dH$ because all elements and subsets of 
$\xi$ belong to $\dV$ by Lemma~\ref{s2w}. 

Conversely if $\xi\in\dH$ is an ordinal in $\dH$ then by 
Lemma~\ref{wo=ss} $\xi$ is equinumerous with a set $w\in\dV.$ Thus 
$w$ admits a well-ordering of length $\xi$ in $\dV$ because subsets 
of sets in $\dV$ belong to $\dV$ by Lemma~\ref{s2w}. Therefore $\xi$ 
is order isomorphic to a set $\xi'\in\dV$ which is an ordinal in 
$\dV.$ This easily implies $\xi=\xi'\in\dV$.\qed\eproof

Let a {\it cardinal\/} mean: an ordinal not equinumerous to a 
\index{cardinal} 
smaller ordinal. 

\bcor
\label{c}
Universes\/ $\dH$ and\/ $\dV$ contain the same cardinals.

The notions of a regular, singular, inaccessible cardinal, and 
the exponentiation of cardinals, are absolute in $\dH$ for $\dV$.
\ecor
\proof If $\kpa\in\dV$ is a cardinal in $\dV$ then at least $\kpa$ 
is an ordinal in $\dH.$ A possible bijection onto a smaller ordinal 
in $\dH$ is effectively coded by a subset of $\kpa\ti\kpa,$ 
therefore it would belong to $\dV.$ The absoluteness holds 
because $\dV$ contains all its subsets in $\dH$ as elements 
by Lemma~\ref{s2w}.\qed

\bdf
\label{doc}
$\Ord$ is the class of all ordinals in $\dH$ (or in $\dV,$ which  
\index{class!ord@$\Ord$ of all ordinals}
is equivalent by the above). Elements of $\Ord$ will be called 
{\it ordinals\/} 
\index{ordinal} 
below. $\Card$ is the class of all cardinals in $\dH$ (or in $\dV$). 
\index{class!card@$\Card$ of all cardinals}
Elements of $\Card$ will be called {\it cardinals\/} below.\qed
\index{cardinal} 
\edf
Ordinals, cardinals in $\dS$ and $\dI$ will be called resp.\ 
\index{ordinal!so@\dd\dS ordinal}%
\index{ordinal!io@\dd\dI ordinal}%
\index{cardinal!sc@\dd\dS cardinal}%
\index{cardinal!ic@\dd\dI cardinal}%
\dd\dS{\it ordinals\/}, \dd\dS{\it cardinals\/}, and 
\dd\dI{\it ordinals\/}, \dd\dI{\it cardinals\/}.

Since $\dV$ models $\ZFC,$ the ordinals satisfy usual theorems, 
e. g. $\Ord$ is well-ordered by the relation: $\al<\ba$ iff 
$\al\in\ba,$ an ordinal is the set of all smaller ordinals, 
$0=\emps$ is the least ordinal, there exist limit ordinals, etc. 
Furthermore the ordinals can be used to define the rank of sets 
in $\dH$ over $\dI,$ the internal subuniverse.

\bdf
\label{irank}
The {\it rank over $\dI,$} $\irk x\in\Ord,$ is defined for each 
\index{rank!irk@$\irk x$}
\index{zzirkx@$\irk x$}
set $x$ in $\dH$ as follows: $\irk x=0$ for internal sets $x,$ 
and $\irk x={\tt sup}_{y\in x}\irk y$ for $x\not\in\dI$.\qed
\edf
(For $O\sq\Ord,$ ${\tt sup}\,O$ is the least ordinal strictly 
bigger than all ordinals in $0$.) This is well defined in $\dH,$ 
by Axiom~\ref{wr} (Weak Regularity).

\subsection{Change of standpoint. Asterisks}
\label{change}

It looks quite natural that $\dS$ and $\dI,$ the classes of all 
resp.\ standard and internal sets, are the principal objects of consideration, e.\ g.\ because $\dS$ is naturally identified 
with the original set universe of ``conventional'' mathematics 
while $\dI$ with an ultrapower of $\dS.$ Following this approach, 
one considers $\dH$ as an auxiliary universe and the notions 
related to $\dH$ as auxiliary notions while the notions related 
to $\dS$ or $\dI$ as primary notions. 

However at this moment it becomes more convenient to treat the 
notions related to $\dH$ and $\dV$ as primary notions,
as in Definition \ref{doc} above. In a sense, $\dV$ is a better 
copy of the ``conventional'' set universe in $\dH$ than $\dS$ 
is, in particular because $\dV,$ unlike $\dS,$ is transitive. 

This change of standpoint leads to an interesting parallel with 
the model theoretic version of nonstandard analysis. 

\bdf
\label{a}
Let, for a set $w\in\dV,$ $\upa w$ denote the set $x\in\dS$ 
\index{zzwast@$\upa w$}
\index{asteriscs}
(unique by Lemma~\ref{s2w}) which satisfies $w=\wx$.\qed
\edf

\bcor
\label{w2s}
$w\map \upa w$ is an\/ \mem-isomorphism\/ $\dV$ onto\/ $\dS$.
 
$\al\in\dV$ is an ordinal (in\/ $\dV$ or in\/ $\dH$) \ iff\/ \ 
$\upa\al$ is an\/ \dd\dS ordinal. 

$\al\in\dV$ is a cardinal (in\/ $\dV$ or in\/ $\dH$) \ iff\/ \ 
$\upa\al$ is an\/ \dd\dS cardinal.\qed
\ecor

\noi
We have approximately the same as what they deal with inside the 
model theoretic approach: $\dH$ corresponds to the basic set 
universe, $\dV$ to a standard model, $\dI$ to its ultrapower (or a 
nonstandard extension of another type); the map $w\map \upa w$ is 
an elementary embedding $\dV$ to $\dI$.

It is an advantage of our treatment that the basic relations in both 
$\dV$ and $\dI$ are of one and the same nature, namely, restrictions 
of the basic relations in $\dH,$ the external universe. 

\subsection{Finite sets and natural numbers}
\label{finsets}

Let a {\it natural number\/} mean: an ordinal smaller than the 
\index{natural number}
\index{set!finite set}
\index{set!nn@$\dN$ of all natural numbers}
least limit ordinal. Let a {\it finite set\/} mean: a set 
equinumerous to a natural number (which is, as usual, equal to 
the set of all smaller numbers). 

The notion of finite set is absolute for $\dV$ in $\dH$ because 
every subset of $\dV$ belongs to $\dV.$ On the other hand, 
$w\in\dV$ is finite in $\dV$ iff $\upa w$ is finite in $\dS,$ by 
Corollary~\ref{w2s}.

\bdf
\label{dn}
$\dN$ is the set of all natural numbers in $\dH$ (or in $\dV,$ 
which is equivalent). Elements of $\dN$ will be called 
{\it natural numbers\/} below. A {\it finite set\/} will mean: \index{natural number}
\index{set!finite set} 
a set finite in the sense of $\dH$ (or $\dV,$ which is 
equivalent provided the set belongs to $\dV.$)\qed
\edf 
Natural numbers in $\dS$ and $\dI$ will be called resp. 
\dd\dS{\it natural numbers\/} (this will become obsolete) and 
\dd\dI{\it natural numbers\/}. 
\index{natural number!snn@\dd\dS natural number}
\index{natural number!inn@\dd\dI natural number}
\index{set!finite set}
\index{set!finite set!fins@\dd\dS finite set}
\index{set!finite set!fini@\dd\dI finite set}
The notions of \dd\dS{\it finite} {\it set\/} and 
\dd\dI{\it finite} {\it set\/} will have similar meaning.


\ble
\label{n}\label{n:s=w}
If\/ $n\in\dN$ then $\upa n=n.$ Therefore the classes\/ 
$\dH,\;\dV,\;\dS$ have the same natural numbers.\qed
\ele
\proof We prove the equality $\upa n=n$ by induction on $n.$ 
Suppose that $\upa n=n$ and prove $\upa(n+1)=n+1.$ Since $n$ 
and $n+1$ are consecutive ordinals, we have $n+1=n\cup\ans{n}$ 
in $\dV.$ We conclude that $\upa(n+1)=n\cup\ans{n}$ in $\dS$ 
by Lemma~\ref{s2w}, in $\dI$ by Transfer, and finally in $\dH$ 
because $\dI$ is transitive in $\dH.$ Thus $\upa(n+1)=n+1$.\qed 

\ble
\label{sf}\label{fin}
Any standard\/ \dd\dS finite set\/ $X$ satisfies\/ $X\sq \dS.$ 
Conversely any finite\/ $X\sq\dS$ is standard and\/ 
\dd\dS finite.  
\ele
\proof Let $X\in\dS$ be \dd\dS finite. Then $X=\ans{f(k):k<n},$ 
where $n$ is an \dd\dS natural number and $f$ is a standard 
function. Then $n=\upa n,$ and $k=\upa k\in\dS$ by Lemma~\ref{n}, 
so every $x=f(k)\in X$ is standard by Transfer. 

To prove the converse, let $X\sq\dS$ be a finite set. Then 
$Y=\ans{\wx:x\in X}$ is a finite subset of $\dV,$ so that 
$Y\in\dV$ by Lemma~\ref{s2w}. The set $\upa Y\in\dS$ is a 
standard \dd\dS finite set, therefore $\upa Y\sq\dS.$ We 
observe that $X$ and $\upa Y$ contain the same standard 
elements, since $\upa(\wx)=x.$ Thus $X=\upa Y$ is a standard 
\dd\dS finite set, as required.\qed


\subsection{Axioms of $\protect\BST$ in the internal subuniverse} 
\label{bst:in:i}

{\it Bounded set theory\/} $\BST$ 
\index{bounded set theory@bounded set theory, $\BST$}
(explicitly introduced by 
Kanovei~\cite{rms}, but very close to the ``internal part'' of a 
theory in~\cite{h78}) is a theory in the \ste-language, 
which includes all of $\ZFC$ (in the \mem-language) together with 
the following axioms:\vspace{2mm} 

\noi{\itt Bounded Idealization} $\aBI$: 
\index{axiom!boun@Bounded Idealization}
$\fstf A\;\exists\,x\in X\;\forall\,a\in A\;\Phi(x,a)\;\;
\llra\;\;\exists\,x\in X\;\fst a\;\Phi(x,a)$;\vspace{2mm}

\noi{\itt Standardization} $\aS$:
\index{axiom!Standardization}
$\;\fst X\;\est Y\;\fst x\;[\,x\in Y\;\;\llra\;\;
x\in X\cj\Phi(x)\,]$;\vspace{2mm}

\noi{\itt Transfer} $\aT$:
\index{axiom!Transfer}
$\;\exists\,x\;\Phi (x)\;\;\lra\;\;\est x\;\Phi(x)$;\vspace{2mm}

\noi{\itt Boundedness\/} $\aB$: 
\index{axiom!Boundedness}
$\;\forall\,x\;\est X\;(x\in X)$.\vspace{2mm}

\noi The formula $\Phi$ must be an \mem-formula in $\aBI$ and 
$\aT,$ and $\Phi$ may contain only standard sets as parameters 
in $\aT,$ but $\Phi$ can be any \ste-formula in $\aS$ and contain 
arbitrary parameters in $\aBI$ and $\aS.$ 
$\fstf A$ means: {\it for all standard finite $A$}. 
\index{zzqfstf@$\fstf$}
$X$ is a standard set in $\aBI$. 

Thus $\aBI$ is weaker than the Idealization $\aI$ of internal 
set theory $\IST$ of Nelson~\cite{ne77} ($\aI$ results by replacing 
in $\aBI$ the set $X$ by the universe of all sets), but the 
Boundedness axiom $\aB$ is added.

We proved in \cite{hyp2} that any model $\dI$ of $\BST$ can be 
enlarged to a $\HST$ model (where it becomes the class of all 
internal sets) by assembling sets along well-founded trees. 
The following is the converse.
\bpro
\label{bstbi}
The class\/ $\dI$ of internal sets in\/ $\dH$ models $\BST$.
\epro
\proof 
\noi Boundedness in $\dI$ follows by the definition of the 
formula $\int.$ The $\BST$ Standardization in $\dI$ follows from 
Axiom~\ref{sza}. Only the $\BST$ Bounded Idealization $\aBI$ 
needs some care.  

Let $\Phi$ be an \mem-formula with internal sets as 
parameters, $X$ a standard set. We have to prove that the 
following is true in $\dI$:\vspace{2mm} 

\noi $\aBI$: \ 
$\fstf A\;\exists\,x\in X\;\forall\,a\in A\;\Phi(x,a)\;\;
\llra\;\;\exists\,x\in X\;\fst a\;\Phi(x,a)$.\vspace{2mm}

\noi
We prove only the direction $\lra;$ the other direction 
follows from the fact that standard \dd\dS finite sets contain 
only standard elements by Lemma~\ref{sf}. 

The main technical problem is to bound the variable $a$ by a 
standard set. We can assume that $\Phi$ contains only one 
parameter, an internal set $p_0,$ which is, by definition, a 
member of a standard set $P.$ For any $a,$ we let 
$Z_a=\ans{\ang{p,x}\in P\ti X:\Phi\upin(x,a,p)}.$ By the $\ZFC$ 
Collection and Transfer, there exists a standard set $A_0$ 
such that $\forall\,a'\,\exists\,a\in A_0\,(Z_a=Z_{a'}).$ We put 
$X_a=\ans{x\in X:\Phi\upin(x,a,p_0)}$ for all $a$.

We verify that the family $\cX$ of all sets $X_a,$ 
${a\in\upSG A_0=\ans{a\in A_0:\st a}},$ satisfies the requirements 
of Axiom~\ref{sat} (Saturation). Indeed each set $X_a$ is  
internal, being 
defined in $\dI$ by an \mem-formula with internal parameters. 
Let $\cX'\sq\cX$ be a finite subset of $\cX.$ By Replacement in 
$\dH,$ there exists a finite set $A\sq A_0$ such that 
$\cX'=\ans{X_a:a\in A}.$ We observe that $A$ is standard and 
\dd\dS finite by Lemma~\ref{fin}. Therefore, by the 
left--hand side of $\aBI,$ the intersection 
$\bigcap\cX'=\bigcap_{a\in A} X_a$ is nonempty, as required. 

Axiom \ref{sat} gives an element $x\in \bigcap\cX.$ We prove 
that $x$ witnesses the right--hand side of $\aBI.$ It suffices to 
check $\fst a'\,\est a\in A_0\,(X_a=X_{a'}).$ Consider a standard 
$a'.$ Then $Z_a=Z_{a'}$ for a standard $a\in A_0$ by the 
choice of $A_0,$ so  
${
X_a=\ans{x:\ang{p_0,x}\in Z_a}=\ans{x:\ang{p_0,x}\in Z_{a'}}=
X_{a'}
}$
.\qed

\subsection{Elementary external sets}
\label{eest}

Let an {\it elementary external set\/} mean a (in $\dI$) \ste-definable  
\index{set!elementary external set}
 subclass 
of an internal set. 
This looks 
unsound, 
but fortunately objects of this type admit a sort of uniform 
description (first discovered in~\cite{iran}), 
given by 
%
\vspace{2mm}

\noi 
$\cC_p=\bigcup_{a\in\upSG A}\,\bigcap_{b\in\upsG B}\,\eta(a,b),$ \ 
\index{zzcp@$\cC_p$}
\begin{minipage}[t]{0.61\textwidth}
where \ $p=\ang{A,B,\eta},$ $A,\,B$ are standard sets, 
$\eta$ being an internal function defined on $A\ti B$.
\end{minipage}
\vspace{2mm}

\noi
If $p\in\dI$ is not of the mentioned form then we set $\cC_p=\emps$.


\bte
\label{tp}
Let\/ $\Phi(x,q)$ be a\/ \ste-formula. The following is a theorem 
of\/ $\BST:$ 
$\forall\,q\;\fst X\;\exists\,p\;(\cC_p=\ans{x\in X:\Phi(x,q)})$.\qed 
\ete
This result (Theorem 16.3 in \cite{iran}, Theorem 2.2 in \cite{hyp2}) 
is an easy consequence of a theorem which asserts that every 
\ste-formula is provably equivalent in $\BST$ to a $\Ss$ 
formula~\footnote
{\ We recall that $\Ss$ denotes the class of all formulas 
\index{zzSi12@$\Ss$}
\index{formula!Sifo@$\Ss$ formula}
$\est a\,\fst b$(\mem-{\itt formula\/}).} 
(Theorem~\ref{bstss} in \cite{hyp1}),  
and a lemma which allows to restrict the two principal quantifiers 
in a $\Ss$ formula by standard sets 
(Lemma~\ref{lem} in \cite{hyp1}, see a corrected proof in 
\cite{hyp3} or a more complicated earlier proof in \cite{iran}, 
Lemma 15.1).

\np

%
\section{Constructibility from internal sets in $\protect\HST$}
\label{str}

Let, as above, $\dH$ be a universe of $\HST,$ $\dS\sq\dI$ be 
the classes of all shandard and internal sets in $\dH,$ $\dV$ 
the condensed subuniverse of $\dH$ introduced in 
Subsection~\ref{conden}. 

The aim of this section is to define an inner class in $\dH$ 
which models $\HST+\neg\;\ip,$ a contribution to the independence 
part of Theorem~\ref{maint}. We shall use $\dL[\dI],$  
{\it the class of all sets constructible from internal sets\/}. 
The main result (Theorem~\ref{li:hst} below) is similar to what 
might be expected in the  
$\ZFC$ case: $\dL[\dI]$ models $\HST$ plus an extra choice--like 
principle. In addition, the isomorphism property $\ip$ 
{\it fails\/} in $\dL[\dI]$. 

$\HST$ is obviously much more cumbersome theory than $\ZFC$ is; 
this leads, in principle, to many additional complications which 
one never meets running the constructibility in $\ZFC$.

On the other hand, it is a certain relief that the initial class 
$\dI$ contains all standard sets. Indeed, since $\dS$ models 
$\ZFC,$ one has, in $\dI,$ an already realized example of the 
constructible hierarchy, essentially of the same length as we are 
looking for in $\dH,$ because $\dS$ and $\dH$ have order isomorphic 
classes of ordinals by Corollary~\ref{w2s}. This allows to use a 
strategy completely different from that in $\ZFC,$ to define 
constructible sets. We introduce $\dL[\dI]$ as the class 
of all sets obtainable in $\dH$ via the procedure of assembling 
sets along wellfounded trees (which starts from sets in $\dI$ and 
involves trees \ste-definable in $\dI$) as they obtain sometimes 
models of fragments of $\ZFC$ from models of 2nd order Peano 
arithmetic.

To make the exposition self-contained, we give a brief review of 
the relevant definitions and results in \cite{hyp3}.

\subsection{Assembling sets along wellfounded trees}
\label{assem}

Let $\Seq$ denote the class of all internal sequences, of 
\index{zzseq@$\Seq$}
arbitrary (but internal) sets, of finite 
length. For $t\in \Seq$ and every set $a,$ 
\index{zztwea@$t\we a$ and $a\we t$}
$t\we a$ is the sequence in $\Seq$ obtained by adjoining $a$ as 
the rightmost additional term to $t.$ The notation $a\we t$ is 
to be understood correspondingly. 

A {\it tree\/} is a non-empty (possibly non--internal) set $T\sq \Seq$ 
\index{tree}
such that, whenever $t',\,t\in \Seq$ satisfy $t'\sq t,$ $t\in T$ 
implies $t'\in T.$ Thus every tree $T$ contains $\La,$ the empty 
\index{zzla@$\La,$ the empty sequence}
sequence, and satisfies $T\sq\dI$. 

$\Max T$ is the set of all \dd\sq{\it maximal\/} in $T$ elements 
\index{zzmaxt@$\Max T$}
$t\in T$.

A tree $T$ is {\it wellfounded\/} ({\it wf tree\/}, in brief) if 
\index{tree!wellfounded, wf tree}
and only if every non-empty (possibly non--internal) set $T'\sq T$ 
contains a \dd\sq maximal element.

\bdf
\label{dfpairs}
\label{dfram}
\label{itf}
Let a {\it wf pair\/} be any pair $\wfp TF$ such that $T$ is a wf 
\index{wf pair}
tree and $F:\Max T\;\lra\;\dI.$ In this case, the 
family of sets $F_T(t),$ $t\in T,$ is defined, using the $\HST$ 
\index{zzftt@$F_T(t)$}
Replacement, as follows:\its
\bit
\item[1)] \ if $t\in \Max T$ then $F_T(t)=F(t)$;\its

\item[2)] \ if $t\in T\setminus\Max T$ then 
$F_T(t)=\ans{F_T(t\we a):t\we a\in T}$.\its
\eit
We finally set $F[T]=F_T(\La)$. 
\index{zzft@$F[T]$}
\qed
\edf
Let, for example, $T=\ans{\La}$ and $F(\La)=x\in\dI.$ Then 
$F[T]=F_T(\La)=x$.


\subsection{Class of elementary external sets}
\label{ees}
      
In particular we shall be interested to study the construction of 
Definition~\ref{dfram} from the point of view of the class 
$\dE=\ans{\cC_p:p\in\dI}$ of all elementary external sets. 
\index{class!e@$\dE$ of all elementary external sets}
(See the definition of $\cC_p$ in Subsection~\ref{eest}.) 


\bpro
\label{e:est}
$\dE$ is a transitive subclass of\/ $\dH.$ $\dE$ models Separation 
in the\/ \ste-language.
\epro
\proof $\dI$ is a model of $\BST,$ see Proposition~\ref{bstbi}. 
It follows that every \ste-definable in $\dI$ subclass of a set 
in $\dI$ has the form $\cC_p$ for some $p\in\dI$ by 
Theorem~\ref{tp}.\qed\eproof

We observe that $\dI\sq\dE,$ and every set $X\in\dE$ satisfies 
$X\sq\dI$.  

\bdf
\label{dfch} $\cH$ is the class of all wf pairs $\wfp TF$ s. t.  
$T,\,F\in \dE$.\qed
\index{class!hc@$\cH$ of all wf pairs}
\edf 
Let $\wfp TF\in\cH.$ Since all sets in $\dE$ are subsets of $\dI,$ 
the set $F[T]$ cannot be a member of $\dE.$ However, 
one can determine, in $\dE,$ different properties which sets of 
the form $F[T]$ ($\wfp TF$ being wf pairs in $\cH$) may have in 
$\dH,$ using the following proposition, proved in \cite{hyp3}.

\bpro
\label{respoc}
$\cH$ is\/ \ste-definable in\/ $\dE$ as a subclass of $\dE\ti\dE.$  
There exist 4--ary\/ \ste-predicates\/ ${}\req{}$ and\/ 
${}\rin{},$ and a binary\/ \ste-predicate\/ $\rst,$ such 
that the following holds for all wf pairs\/ $\wfp TF\,,\;\wfp RG$ 
in\/ $\cH\;:$
\dm
\bay{cccc}
F[T]=G[R]&\hbox{ iff }&\hbox{it is true in\/ }\;\dE\;
\hbox{ that} & \wfp TF\req \wfp RG\,;\\[2mm]

F[T]\in G[R]&\hbox{ iff }&\hbox{it is true in\/ }\;\dE\;
\hbox{ that} & \wfp TF\rin \wfp RG\,;\\[2mm]

\st F[T] &\hbox{ iff }& \hbox{it is true in\/ }\;\dE\;
\hbox{ that} & \rst \wfp TF\,.
\eay
\dm
\epro
\proof (An outline. See a complete proof in 
\cite{hyp3}, Section 3.) 
To prove the definability of $\cH$ in $\dE,$ 
it suffices to check that if $T\in\dE$ is a wf tree in 
$\dE$ then\/ $T$ is a wf tree in the sense of $\dH,$ too. 
Since $\dE$ models Separation, the wellfoundedness of $T$ in $\dE$ 
allows to define, in $\dE,$ the rank function $\rho$ from $T\,$into 
\dd\dS ordinals. Such a function proves that $T$ is 
wellfounded in $\dH$.

The formula $\wfp TF\req \wfp RG$ expresses the 
existence of a computation of the truth values of 
equalities $F_T(t)=G_R(r),$ where $r\in R$ and $t\in T,$ 
which results in \ {\sf true} \ for the equality 
$F_T(\La)=G_R(\La).$ The other two predicates are simple 
derivates of ${}\req{}$.\qed

\subsection{The sets constructible from internal sets}
\label{main}

The following definition introduces the sets constructible 
from internal sets. 

\bdf
\label{defli}
$\dL[\dI]=\ans{F[T]:\wfp TF\in \cH}$.\qed
\index{class!li@$\dL[\dI]$ of all sets constructible from 
internal sets}
\edf
In principle this does not look like a definition of 
constructibility.~\footnote
{\ This approach to the constructibility from internal sets in 
$\HST$ was introduced in~\cite{hyp3}. For any infinite cardinal 
$\kpa$ we defined the class $\dI_\kpa\sq\dI$ of all internal sets 
which are elements of standard sets of cardinality $\<\upa\kpa$ in 
$\dS,$ and then the class $\dH_\kpa$ (could be denoted by 
$\dL[\dI_\kpa]$) of all sets constructible in this sense from sets 
in $\dI_\kpa.$ We proved in \cite{hyp3} that $\dH_\kpa$ models a 
certain \dd\kpa version of $\HST;$ the proof of Theorem~\ref{li:hst} 
below in part copies some arguments from \cite{hyp3}. But we did 
not prove results similar to statements \ref{gdc} and \ref{II} 
in \cite{hyp3}.} 
However it occurs that $\dL[\dI]$ is the least class in $\dH$ 
which contains all internal sets and satisfies $\HST,$ a sort of 
characterization of what in general the class $\dL[\dI]$ should be.    
Anyway this gives the same result as the ordinary definition of 
constructible sets, but with much less effort in this particular 
case.

To formulate the theorem, let us recall some notation related 
to ordered sets. A subset $Q\sq P$ of a p.\ o. set $P$ is called 
{\it open dense\/} in $P$ iff 
\index{set!open dense}
$1)$$\forall\,p\in P\;\exists\,q\in Q\;(q\<p)$ and 
$2)$$\forall\,p\in P\;\forall\,q\in Q\;(p\<q\;\lra\;p\in Q)$.

\bdf
\label{compl}
Let $\kpa$ be a cardinal. A p.\ o.\ set $P$ is \dd\kpa 
{\it closed\/} iff  
\index{set!kclosed@\dd\kpa closed}
\index{set!kdistributive@\dd\kpa distributive}
every decreasing chain $\ang{p_\al:\al<\kpa}$ (i.\ e. 
$p_\al\<p_\ba$ whenever $\ba<\al<\kpa$) in $P$ has a lower bound 
in $P.$ A p.\ o.\ set $P$ is \dd\kpa{\it distributive\/} iff an 
intersection of \dd\kpa many open dense subsets of $P$ is dense in 
$P$.\qed
\edf
The distributivity is used in the practice of forcing as a 
condition which prevents new subsets of sets of certain cardinality 
to appear in generic extensions. We shall use it with the aim to 
preserve Standardization in the extensions. In $\ZFC$ a closed set 
is distrubutive, of course, but this simple fact is based on Choice. 
$\HST$ does not include a sufficient amount of Choice, but we can 
prove the implication to be true in $\dL[\dI]$. 

\bte
\label{li:hst}
$\dL[\dI]$ is a transitive subclass of\/ $\dH,$ $\dE\sq\dL[\dI],$ 
$\dL[\dI]$ models\/ $\HST$ and models the following two additional 
postulates$:$\its
\ben
\def\theenumi{P--\arabic{enumi}}
\def\labelenumi{{\rmt\theenumi}.}
\itla{gdc} 
For each cardinal\/ $\kpa,$ every\/ \dd\kpa closed p.\ o. set 
is \dd\kpa distributive.\its

\itla{II}  
The isomorphism property\/ $\ip$ {\rm (see Introduction)} fails. 
\een
\ete
Take notice that any model $H\sq\dH$ of $\HST$ transitive in $\dH$ 
and containing all internal sets satisfies $\dV\sq H$ and has 
the same classes of ordinals and cardinals as $\dH$ and $\dV$ do, 
by Corollary~\ref{o}. Therefore the cardinals and preassumed 
ordinals in 
item~\ref{gdc} are those given by Definition~\ref{doc} in $\dH$.

Furthermore, since $\dI$ is transitive in $\dH,$ the construction 
of $\dL[\dI]$ is absolute for any such class $H.$ In other 
words, $\dL[\dI]$ is actually the least transitive class in $\dH$ 
which contains all internal sets and models $\HST$.  
\eproof

\proof To prove the transitivity, let $X=F[T]\in\dL[\dI],$ where
$\wfp TH\in\cH.$ If $T=\ans{\La}$ then $X=F(\La)\in\dI$ by 
definition, hence all elements of $X$ are internal. Suppose that 
$T\not=\ans{\La}.$ Then $\Min T=\ans{a:\ang{a}\in T}$ is 
\index{zzmint@$\Min T$}
a non-empty set. For $a\in\Min T,$ we put $T^a=\ans{t:a\we t\in T}$ 
and $F^a(t)=F(a\we t)$ for all $a\we t\in\Max T.$ Obviously 
$\wfp{T^a}{F^a}$ is a wf pair for all $a;$ moreover $T^a$ and 
$F^a$ belong to $\dE$ by Proposition~\ref{e:est}, so that in fact 
$\wfp{T^a}{F^a}\in\cH$ and $x_a=T^a[F^a]\in\dL[\dI].$ On the other 
hand, $X=\ans{x_a:a\in\Min T}$.

To prove $\dE\sq\dL[\dI],$ let $A\in\dE.$ We define 
$T=\ans{\La}\cup\ans{\ang{a}:a\in A}$ and set $F(\ang{a})=a$ for 
all $a\in A.$ Then $\wfp TF$ is a wf pair; furthermore 
$\wfp TF\in\cH$ by Proposition~\ref{e:est}, so that 
$A=F[T]\in\dL[\dI],$ as required.\eproof

Let us prove three auxiliary claims which will be used below.\eproof

\noi{\bf Claim 1} \ {\it In\/ $\dL[\dI],$ every set is a functional 
image of a standard set\/}.\eproof

\proof
Suppose that $\wfp TF\in\cH,$ 
$X=F[T]\in \dL[\dI].$ Then $A=\Min T$ contains only internal 
elements. By Lemma~\ref{hstb}, there exists a standard set $U$ 
such that $A\sq U.$ For $a\in A,$ we define $f(a)=x_a.$ 
For $a\in U\setminus A,$ let $f(a)=f(a_0),$ where $a_0$ is a 
fixed element of $A.$ Then $f$ maps $U$ onto $X.$  
Proposition~\ref{e:est} allows to transform the given wf pair 
$\wfp TF$ to a wf pair $\wfp RG\in\cH$ such that $f=G[R],$ which 
proves $f\in\dL[\dI],$ as required.\bult\eproof 

\noi{\bf Claim 2} \ {\it Every set\/ $X\sq\dL[\dI]$ of standard 
size in\/ $\dH$ belongs to $\dL[\dI]$.}\eproof

\proof
Lemma~\ref{exten} gives an internal function $f$ defined on 
a standard set $A$ so that, for every 
$a\in\upSG A,$ $f(a)=\ang{\tau_a,\sg_a}$ is an internal pair, the 
sets $T_a=\cC_{\tau_a}$ and $F_a=\cC_{\sg_a}$ satisfy 
$\wfp{T_a}{F_a}\in\cH,$ and $X=\ans{F_a[T_a]:a\in\upSG A}$.   

Let $T=\ans{\La}\cup\ans{a\we t:a\in\upSG A\cj t\in T_a}$ and 
$F(a\we t)=F_a(t)$ for all $a\in\upSG A$ and $t\in\Max T_a.$ Thus 
both $T$ and $F$ are \mem-definable in $\dE$ using only the 
internal $f$ as a parameter. Therefore $F$ and $T$ belong to $\dE,$ 
by Proposition~\ref{e:est}. On the other hand, by definition 
$\wfp TF$ is a wf pair -- thus $\wfp TF\in \cH$ -- and 
$F[T]=\ans{F_a[T_a]:a\in\upSG A}=X\in\dL[\dI]$.\bult\eproof

\noi{\bf Claim 3} \ {\it Every set\/ $Z\in\dL[\dI],$ $Z\sq\dI$ 
belongs to $\dE$}.\eproof

\proof
Suppose that $Z=F[T],$ where $\wfp TF\in\cH,$ in particular 
$T,\,F\in\dE.$ We set $T^x=\ans{\La}$ 
and $F^x(\La)=x$ for all internal $x;$ then obviously 
$\wfp{T^x}{F^x}\in\cH$ and $F^x[T^x]=x.$ Then  
$Z=\ans{x\in\dI:F^x[T^x]\in F[T]}.$ Using propositions \ref{respoc} 
and \ref{e:est}, we finally obtain $Z\in\dE$.\bult

\subsubsection*{Isomorphism property}
 
We prove statement~\ref{II} of the theorem. Since any two infinite 
sets are elementarily equivalent as 
structures of the language containing nothing except the equality, 
the following lemma implies the negation of $\ip$.\eproof

\noi{\bf Lemma} \ {\it Any two internal\/ \dd\dI infinite sets 
of different\/ \dd\dI cardinalities are non--equinumerous in $\dL[\dI]$.}\eproof

\proof Suppose that  $\card X<\card Y$ in $\dI,$ and 
$f\in\dL[\dI]$ maps $X$ into $Y.$ We have to prove that $\ran f$ is 
a proper subset of $Y.$ The following argument is a part of a 
more complicated reasoning in Kanovei~\cite{jsl}. 
(where the $\IST$ case was considered). 

It follows from Claim 3 that $f=\cC_p$ for some 
$p\in\dI,$ so that there exist standard sets $A,\,B$ and an 
internal set $W\sq X\ti Y\ti A\ti B$ such that  
\dm 
f(x)=y 
\hspace{7mm}\hbox{iff}\hspace{7mm}
\est a\in A\;\fst b\in B\;W(x,y,a,b),
\dm 
for all $x\in X,\;y\in Y.$ Since $\dI$ models $\BST,$ there exists 
an \dd\dI finite (but perhaps infinite) set $Z\in\dI$ containing 
all standard elements of $A$ and~$B$. 

We put $F(x,y)=\ans{\ang{a,b}\in Z\ti Z: W(x,y,a,b)}$ 
for $x\in X,\;y\in Y.$ Then obviously $f(x)=y$ iff $f(x')=y',$ 
provided $x,\,x'\in X$ and $y,\,y'\in Y$ satisfy 
$F(x,y)=F(x',y')$.

On the other hand $F$ is an internal function, taking values in 
an \dd\dI finite set $\cP(Z\ti Z).$ Therefore arguing in $\dI$ as 
an \mem-model of $\ZFC$ we obtain a set $Y'$ of \dd\dI cardinality 
$\card Y'\<\card X$ (here the \dd\dI infinity of $X$ is used) such 
that for each $x\in Z$ and $y$ there exists $y'\in Y'$ satisfying 
$F(x,y)=F(x,y').$ In other words $\ran f\sq Y'.$ Finally 
$\card Y'<\card Y$ in $\dI,$ so that $Y'$ is a proper subset of 
$Y,$ as required.\bult\eproof

{\it Remark\/}. Internal {\it finite\/} sets of different 
\dd\dI cardinality in $\dI$ can become equinumerous in $\dL[\dI].$ 
For instance if $n$ is a nonstandard \dd\dI natural number then 
clearly sets containing $n$ and $n+1$ elements are equinumerous 
in an external universe. It follows from a result of 
Keisler, Kunen, Miller, and Leth~\cite{kkml} that if 
$n\<k\<sn$ in $\dI,$ for a standard natural $s,$ then  
sets containing $n$ and $k$ elements are equinumerous 
in an external universe.

\subsubsection*{Verification of $\protect\HST$ axioms} 

We verify the $\HST$ axioms  
in $\dL[\dI].$ 
Axioms of group \ref{it1} in Subsection~\ref{hst:ax}, 
Standardization, Saturation, and Weak Regularity are simply 
inherited from $\dH$ because $\dI\sq\dE\sq\dL[\dI]$. 

To prove standard size Choice in $\dL[\dI],$ we first get a 
choice function in $\dH.$ The function is a standard size subset 
of $\dL[\dI],$ therefore belongs to $\dL[\dI]$ by Claim 2. This 
argument also verifies Dependent Choice in $\dL[\dI]$.

Thus actually only the $\ZFC$ axioms included in group~\ref{hzfc} 
need a consideration. 
Among them, only Separation and Collection do need a serious 
verification; the rest of the axioms are either inherited from 
$\dH$ or proved by elementary transformations of wf pairs involved. 

Let us check {\it Separation\/} in $\dL[\dI].$ Let $X=F[T],$ 
$\wfp TF\in\cH,$ and all parameters in a \ste-formula $\Phi(x)$ 
belong to $\dL[\dI].$ We have to prove that the set 
$X'=\ans{x\in X:\dL[\dI]\models\Phi(x)}$ also belongs to $\dL[\dI]$. 

Suppose that $T\not=\ans{\La}.$ Then the set 
$A=\Min T=\ans{a:\ang{a}\in T}$ is non-empty. For $a\in A,$ we 
define a wf pair $\wfp{T^a}{F^a}\in\cH$ as in the proof of 
transitivity above. Then 
$X=\ans{F^a[T^a]:a\in A}.$ Let $A'$ be the set of all $a\in A$ 
such that $\Phi(F^a[T^a])$ is true in $\dL[\dI].$ Then $A'$ is 
\ste-definable in $\dE$ by Proposition~\ref{respoc}%
, therefore $A'\in\dE$ by Proposition~\ref{e:est}. 
Let us define $T'=\ans{\La}\cup\ans{a\we t\in T:a\in A'},$ and 
$F'(a\we t)=F(a\we t)$ for each $a\we t\in\Max T'.$ Then 
$\wfp{T'}{F'}\in\cH$ by Proposition~\ref{e:est}, 
and $X'=F'[T']$.
 
Suppose that ${T=\ans{\La}},$ so that $X=F(\La)$ is internal. Then 
$X'$ is definable in $\dE$ by Proposition~\ref{respoc}, therefore 
$X'\in\dE$ by Proposition~\ref{e:est}. This implies 
$X'\in\dL[\dI],$ since $\dE\sq\dL[\dI],$ see above. 

Let us check {\it Collection\/} in $\dL[\dI].$ By the $\HST$ 
Collection in $\dH,$ it suffices to verify the following: if a set 
$X\in\dH$ satisfies $X\sq\dL[\dI]$ then 
there exists $X'\in\dL[\dI]$ such that $X\sq X'.$ Using Collection 
again and definition of $\dL[\dI],$ 
we conclude that there exists a set $P\sq\dI$ such that 
\dm
\forall\,x\in X\;\exists\,a=\ang{p,q}\in P\;(\wfp{\cC_p}{\cC_q}\in 
\cH\cj x=x_a=\cC_q[\cC_p])\,.
\dm
By Lemma~\ref{hstb}, there exists a standard set $S$ such that 
$P\sq S.$ Let us define $A$ to be the set of all pairs 
$a=\ang{p,q}\in P$ such that $\wfp{\cC_p}{\cC_q}\in \cH.$ Then 
$A\in\dE,$ as above. 

We set $T=\ans{\La}\cup\ans{a\we t:a=\ang{p,q}\in A\cj t\in\cC_p}$ 
and $F(a\we t)=\cC_q(t)$ whenever $\ang{p,q}\in S$ and 
$t\in\Max\cC_p.$ Then both $T$ and $F$ belong to $\dE$ by 
Proposition~\ref{e:est}. Since obviously $\wfp TF$ is a wf pair, 
we conclude that $\wfp TF\in\cH$ and $X'=F[T]\in\dL[\dI].$ On the 
other hand, $X\sq X'$.

\subsubsection*{Distributivity} 

Let $\kpa$ be a cardinal, 
$P=\ang{S;\<}\in\dL[\dI]$ be a \dd\kpa closed in $\dL[\dI]$ 
partial order, on a {\it standard\/} (Claim 1) set $S.$
Consider a family  
$\ang{D_\al:\al<\kpa}\in\dL[\dI]$  
of open dense subsets of $\ang{S;\<}.$ Let us prove 
that the intersection $\bigcap_{\al<\kpa}D_\al$ is dense in $P.$ 
Let $\ux\in S.$ We have to find an element  
$x\in S,$ $x\<\ux$ such that $x\in\bigcap_{\al<\kpa} D_\al$. 

To simplify the task, let us first correct the order. 
For $x\in S,$ let $\al(x)$ denote the largest 
$\al\<\kpa$ such that $x\in D_\ba$ for all $\ba<\al.$ For 
$x,\,y\in S,$ we let $x\prec y$ mean: $x\<y,$ and either 
$\al(y)<\al(x)<\kpa$ or $\al(x)=\al(y)=\kpa$. 

Now it suffices to obtain a \dd\prec decreasing \dd\kpa sequence 
$\fx=\ang{x_\al:\al<\kpa}$ of 
$x_\al\in S,$ satisfying $x_0\<\ux.$ 
Indeed, then $x_\al\in D_\ba$ for all $\ba<\al<\kpa.$ Furthermore 
$\fx\in\dL[\dI]$ by Claim 2. 
It follows that some $x\in S$ is $\<$ each $x_\al$ because $P$ 
is \dd\kpa closed in $\dL[\dI].$ Then 
$x\in\bigcap_{\al<\kpa} D_\al$. 

The order relation $\prec$ belongs to $\dL[\dI]$ 
by the already verified Separation in $\dL[\dI].$ It follows that 
$\prec$ is $\cC_p$ for some internal $p,$ by Claim 3; in other 
words, there exist standard sets $A',\,B'$ and an internal set 
$Q\sq A'\ti B'\ti S^2$ such that  
$x\prec y$ iff ${\est a\in A'\;\fst b\in B'\;Q(a,b,x,y)}$ --- 
for all $x,\,y\in S$. 

We have $A'=\upa A$ and $B'=\upa B,$ for some $A,\,B\in\dV,$ by 
Corollary~\ref{w2s}. Let 
$Q_{ab}=\ans{\ang{x,y}\in S^2:Q(\upa a,\upa b,x,y)}$ for $a\in A,$ 
$b\in B.$ Then, in $\dH,$ $x\prec y$ iff 
$\exists\,a\in A\;\forall\,b\in B\;Q_{ab}(x,y),$ and $Q_{ab}$ are 
internal sets for all~$a,\,b$.

The principal idea of the following reasoning can be traced down 
to the proof of a choice theorem in Nelson~\cite{ne88}: we divide 
the problem into a choice argument in the \mem-setting and a 
saturation argument. 

Let us say that ${a\in A}$ {\it witnesses\/} ${x\prec y}$ 
iff we have ${\forall\,b\in B\,Q_{ab}(x,y)}$. 

For any $\al\<\kpa,$ we let $\cA_\al$ be the family of all 
functions $\fa:\al\ti\al\,\lra\,A$ such that there exists a 
function $\fx:\al\,\lra\,S$ satisfying $\fx(0)\<\ux$ and the 
requirement that $\fa(\da,\ga)\in A$ witnesses 
$\fx(\ga)\prec \fx(\da)$ whenever $\da<\ga<\al$. 

We observe that, by Lemma~\ref{s2w}, each function 
$\fa\in\bigcup_{\al\<\kpa}\cA_\al,$ every set $\cA_\al,$ and the 
sequence $\ang{\cA_\al:\al\<\kpa}$ belong to $\dV$. 

It suffices to prove that $\cA_\kpa\not=\emps$. 

Since the sequence of sets $\cA_\al\;\;(\al\<\kpa)$ belongs to 
$\dV,$ a $\ZFC$ universe, the following facts 1 and 2 immediately 
prove $\cA_\kpa\not=\emps$.\eproof

\noi
{\bf Fact\ 1} \ {\it If\/ $\al<\kpa$ and\/ $\fa\in\cA_\al$ then 
there exists\/ $\fa'\in\cA_{\al+1}$ extending $\fa$}.\eproof

\proof By definition there exists an decreasing \dd\al chain 
$\fx:\al\,\lra\,S$ such that $\fa(\da,\ga)$ witnesses 
$\fx(\ga)\prec \fx(\da)$ whenever $\da<\ga<\al.$ Since $P$ is 
\dd\kpa closed, some $x\in S$ is  
$\<\fx(\da)$ for each $\da<\al.$ By the density of the sets 
$D_{\ba},$ we can assume that in fact $x\prec \fx(\da)$ for all 
$\da<\al.$ Using the standard size Choice 
in $\dH,$ we obtain a function $f:\al\,\lra\,A$ such that $f(\da)$ 
witnesses $x\prec \fx(\da)$ for each $\da<\al.$ We define 
$\fa'\in\cA_{\al+1}$ by\linebreak[3] 
$\fa'(\da,\ga)=\fa(\da,\ga)$ whenever 
$\da<\ga<\al,$ and $\fa'(\da,\al)=f(\da)$ for $\da<\al$.\bult\eproof

\noi
{\bf Fact\ 2} \ {\it If\/ $\al<\kpa$ is a limit ordinal and a 
function\/ $\fa:\al\ti\al\,\lra\,A$ satisfies\/ 
$\fa\res (\ba\ti\ba)\in\cA_\ba$ for all\/ $\ba<\al$ then
$\fa\in\cA_\al$}.\eproof

\proof Suppose that $\da<\ga<\al$ and $b\in B.$ We let 
$\Xi_{b\da\ga}$ be the set of all internal functions 
$\xi:\upa\al\,\lra\,S$ such that 
$Q_{\fa(\da,\ga)\,b}(\xi(\upa\da),\xi(\upa\ga)).$ 
The sets $\Xi_{b\da\ga}$ are internal because so are all $Q_{ab}$. 

We assert that the intersection
$
\Xi_\ba=
{\textstyle\bigcap_{\;b\in B;\;\,\da<\ga<\ba\;}}\Xi_{b\da\ga}
$
is non-empty, for any $\ba<\al.$ Indeed, since 
$\fa\res (\ba\ti\ba)\in\cA_\ba,$ there exists a function 
$\fx:\ba\,\lra\,S$ such that 
$\fa(\da,\ga)$ witnesses $\fx(\ga)\prec \fx(\da)$ whenever 
$\da<\ga<\ba.$ By the Extension lemma (Lemma~\ref{exten}) there 
exists an internal function $\xi,$ defined on $\upa\al$ and 
satisfying $\xi(\upa\ga)=\fx(\ga)$ for all $\ga<\al.$ Then 
$\xi\in\Xi_\ba$. 

Then the total intersection 
$\Xi=\bigcap_{\;b\in B;\;\,\da<\ga<\al\;}\Xi_{b\da\ga}$ is 
non-empty by Saturation in $\dH.$ Let $\xi\in \Xi.$ By 
definition, we have 
$Q_{\fa(\da,\ga)\,b}(\xi(\upa\da),\xi(\upa\ga))$ whenever 
$\da<\ga<\al$ and $b\in B.$ Let 
$\fx(\da)=\xi(\upa\da)$ for all $\da<\al.$ Then 
$Q_{\fa(\da,\ga)\,b}(\fx(\da),\fx(\ga))$ holds whenever 
$\da<\ga<\al$ and $b\in B.$ In other words, $\fx$ shows that 
$\fa\in\cA_\al,$ as required.\bult\qed

\np

%
\section{Forcing over models of $\protect\HST$}
\label{f}

The proof of the consistency part of Theorem~\ref{maint} involves 
forcing. This section shows how in general one can develop forcing 
for $\HST$ models.

It is a serious problem that the membership relation is not 
well-founded in $\HST.$ This does not allow to run forcing over a 
$\HST$ model entirely in the $\ZFC$ manner: for instance the 
induction on the \mem-rank, used to define the forcing relation 
for atomic formulas, does not work. 

However this problem can be solved, using the axiom of Weak 
Regularity, or well-foundedness over the internal subuniverse 
$\dI$ (\ref{wr} in Subsection~\ref{hst:ax}). 

We shall assume the following.\its

\ben
\def\theenumi{$(\fnsymbol{enumi})$}
\def\labelenumi{\theenumi}
\addtocounter{enumi}{1}
\itla{dhis}\msur
$\dH$ is a model of $\HST$ in a wider set  universe~\footnote
{\ Say a $\ZFC$ universe. The membership relation $\inh$ in $\dH$ 
may have nothing in common with the true membership in the wider 
universe.} 
. $\dS\sq\dI$ and $\dV$ are resp. the classes of all standard and 
\index{class!s@$\dS$ of all standard sets}%
\index{class!i@$\dI$ of all internal sets}%
\index{class!h@$\dH$ of all external sets}%
\index{class!v@$\dV$ of all well-founded sets}%
internal sets in $\dH,$ and the condensed subuniverse defined 
in $\dH$ as in Subsection~\ref{conden}.\its 

\itla{dhwf}\msur 
$\dH$ {\it is well-founded over}~$\dI$ in the sense that the 
\index{class!hwi@$\dH,$ well-founded over $\dI$}%
ordinals of $\dH$ ( = those of $\dV$) are well-founded in the wider 
universe. (Or, equivalently, $\dV$ is a well-founded \mem-model.)
\its
\een
In this case, we shall study generic extensions of $\dH$ viewing 
$\dH$ as a sort of \dd\ZFC like model 
with urelements; internal sets playing the role of urelements. 
Of course internal sets behave not completely like urelements; 
in particular they participate in the common membership relation. 
But at least this gives an idea how to develop forcing in this 
case: the extension cannot introduce new internal sets 
(therefore neither new standard sets nor new well-founded sets --- 
members 
of the condensed subuniverse $\dV.$) Thus, in the frameworks of 
this approach, we can expect to get only new non--internal sets.   

One more problem is the Standardization axiom. Since new standard 
sets cannot appear, a set of standard size, in particular a set in 
$\dV,$ cannot acquire new {\it subsets\/} in the extension. To obey 
this restriction, we apply a classical forcing argument: if the 
forcing notion is standard size distributive then no new standard  
size subsets of $\dH$ appear in the extension. 

\subsection{The extension} 
\label{gen}

Let $\dP=\ang{\dP;\<}$ be a partially ordered set in $\dH$ --- the 
{\it forcing notion\/}, containing the maximal element $\onep.$ 
Elements of $\dP$ will be called {\it (forcing) conditions\/} and 
denoted, as a rule, by letters $p,\,q,\,r$.\nopagebreak 

The inequality $p\<q$ means that $p$ is a {\it stronger\/} condition.
\index{forcing!anot@notion}
\index{forcing!acond@condition}
\index{forcing!acond@condition!stronger}

Let $\brx=\ang{0,x}$ for any set $x\in\dH.$ $\brx$ will be the 
\index{zzxbr@$\brx$}
``name'' for $x.$ We define $\cN_0=\ans{\brx:x\in\dH}.$  
For $\al>0,$ we let 
$\cN_\al=\ans{a:a\sq\dP\ti\bigcup_{\ba<\al}\cN_\ba}.$ 

We observe that ``names'' in $\cN_0$ never appear again at 
higher levels. 

We define, in $\dH,$ $\cN=\cN[\dP]=\bigcup_{\al\in\Ord}\cN_\al,$ 
the class of \dd\dP``names'' 
\index{class!n@$\cN$ of ``names''}
for elements in the planned extension $\dH[G].$ (We recall 
that the class $\Ord$ of all ordinals in $\dH$ was introduced in 
Subsection~\ref{orcar}. It follows from \ref{dhwf} that \dd\dH 
ordinals can be identified with an initial segment of the true 
ordinals in the wider universe.)
For $a\in\cN,$ we let $\nrk a$ (the {\it name--rank\/} of $a$) 
\index{rank!nrk@$\nrk x$}
\index{zznrkx@$\nrk x$}
denote the least ordinal $\al$ such that $a\in\cN_\al$. 

Suppose that $G\sq \dP$ (perhaps $G\not\in\dH$). We define a set 
$a[G]$ in the wider universe 
for each ``name'' $a\in\cN$ by induction on $\nrk a$ as follows. 

First of all, we put $a[G]=x$ in the case when $a=\brx\in\cN_0$. 

Suppose that $\nrk a>0.$ Following the $\ZFC$ approach, we would 
define 
\dm
a[G]=\ans{b[G]:\exists\,p\in G\;(\ang{p,b}\in a)}\,.\eqno{(\ast)}
\dm
However we face a problem: a set $a[G]$ defined this way may 
contain the same elements as some $x\in\dH$ \dd\inh contains in 
$\dH,$ so that $a[G]$ and $x$ must be somehow identified in 
$\dH[G]$ in order not to conflict with Extensionality. This problem 
is settled as follows. We define, as above, for 
$a\in\cN\setminus\cN_0$, 
\dm
a'[G]=\ans{b[G]:\exists\,p\in G\;(\ang{p,b}\in a)}.
\index{zzagp@$a'[G]$}
\index{zzag@$a[G]$}
\dm
If there exists $x\in\dH$ such that  
$y\in a'[G]$ iff $y\in \dH\cj y\inh x$ for each $y,$ 
then we let $a[G]=x.$ Otherwise we put $a[G]=a'[G]$. 
(Take notice that if $\ang{p,b}\in a\in\cN$ for some $p$ then 
$\nrk b<\nrk a,$ so that $a[G]$ is well defined for all 
$a\in\cN,$ because $\dH$ is 
assumed to be well-founded over $\dI$.)

We finally set $\dH[G]=\ans{a[G]:a\in\cN}$. 
\index{zzhg@$\dH[G]$}
\index{class!hg@$\dH[G],$ generic extension}

We define the {\it membership\/} $\ing$ in $\dH[G]$ as follows: 
\index{membership relation!ing@$\ing$ on $\dH[G]$}
\index{membership relation!inh@$\inh$ on $\dH$}
$x\ing y$ iff either $x,\,y$ belong to $\dH$ and $x\inh y$ in 
$\dH,$ or $y\not\in\dH$ and $x\in y$ in the sense of 
the wider universe. We define the {\it standardness\/} in $\dH[G]$ 
by: $\st x$ iff $x\in\dH$ and $x$ is standard in $\dH$. 

\bdf
\label{plain}
A \ste-structure $\dH'$ is a {\it plain extension\/} of $\dH$ 
\index{plain extension}
iff $\dH\sq\dH',$ $\dH$ is an \dd\inhp transitive part of $\dH',$ 
${\inh}={{\inhp}\res\dH},$ and the standard (then also internal) 
elements in $\dH'$ and $\dH$ are the same.\qed
\edf
It is perhaps not true that $\dH[G]$ models $\HST$ independently 
of the choice of the notion of forcing $\dP.$ To guarantee 
Standardization in the extension, new subsets of standard size 
``old'' sets cannot appear. Standard size distributivity 
provides a sufficient condition.

\bdf
\label{ssdis}
A p.\ o.\ set $\dP$ is {\it standard size closed\/} iff it is 
\index{set!standard size closed}
\dd\kpa closed for every cardinal $\kpa.$ A p.\ o.\ set $\dP$ is 
{\it standard size distributive\/} iff it is 
\index{set!standard size distributive}
\dd\kpa distributive for every cardinal $\kpa$.\qed 
\edf

\bte
\label{is}
\label{t:hg}
Let, in the assumptions\/ \ref{dhis} and\/ \ref{dhwf}, $\dP\in\dH$ 
be a p.\ o. set and $G\sq\dP$ be\/ \dd\dP generic over $\dH.$ 
Then\/ $\dH[G]$ is a plain extension of\/ $\dH$ containing\/ $G$ 
and satisfying Extensionality. If in addition the notion of 
forcing\/ $\dP$ is standard size distributive in\/ $\dH$ then\/ 
$\dH[G]$ models $\HST$.
\ete
\proof $\dH\sq \dH[G]$ because $\brx[G]=x$ by definition. 
Furthermore putting $\brG=\ans{\ang{p,\brp}:p\in \dP},$ we get 
$\brG[G]=G$ for all $G\sq\dP,$ so that $G$ also belongs to $\dH[G].$ 
The membership in $\dH$ is the restriction of the one in $\dH[G]$ 
by definition, as well as the \dd\ing transitivity of $\dH$ in 
$\dH[G]$ and the fact that the standard sets are the same in 
$\dH$ and $\dH[G]$.

To prove Extensionality, let $a[G],\,b[G]\in\dH[G]$ 
\dd\ing contain the same elements in $\dH[G];$ we have to prove 
that $a[G]=b[G].$ If $a[G]=A\in\dH$ then $a[G]$ \dd\ing contains 
the same elements in $\dH[G]$ as $A$ \dd\inh contains in $\dH,$ so 
that $b[G]=A$ by definition. The case $b[G]\in\dH$ is similar. 
If $a[G]\not\in\dH$ and $b[G]\not\in\dH$ then by definition $a[G]=a'[G]=b'[G]=b[G]$.

To proceed with the proof of the theorem, we have to define forcing. 

\subsection{The forcing relation} 
\label{fo}

We argue in the model $\dH$ of $\HST$ in this subsection.
 
Let $\dP\in\dH$ be a p.\ o. set. 
\index{forcing!arel@relation}
The aim is to define the forcing relation ${\fo}={\fo}_\dP,$ used 
\index{forcing!fo@$\fo$}
as $p\fo\Phi,$ where $p\in \dP$ while $\Phi$ is a \ste-formula with 
``names'' in $\cN$ as parameters. 

First of all let us consider the case when $\Phi$ is an atomic 
formula, $b=a$ or $b\in a,$ where $a,\,b\in\cN.$ The definition 
contains several items.
\ben
\def\theenumi{F-\arabic{enumi}}
\itla{f1}
We define: \ $p\fo \brx=\bry$ \ iff \ $x=y,$ \ and \ 
$p\fo \brx\in \bry$ \ iff \ $x\in y$.
\een
Let $a,\,b\in\cN.$ We introduce the auxuliary relation
\dm
p\sfo b\in a\hspace{4mm}\hbox{iff}\hspace{4mm}\left\{
\index{forcing!sfo@$\sfo$}
\bay{cl}
\exists\,y\in x\;(b=\bry) & \hbox{whenever }\,a=\brx\in\cN_0\\[2mm]

\exists\,q\>p\;(\ang{q,b}\in a) & \hbox{otherwise}
\eay
\right.
\dm
Note that $p\sfo b\in a$ implies that either $a,\,b\in\cN_0$ or 
$\nrk b<\nrk a$.
\ben
\def\theenumi{F-\arabic{enumi}}
\addtocounter{enumi}{1}
\itla{f=}\msur
$p\fo a=b$ \ iff \ for every condition $q\<p$ the following 
holds:\vspace{1mm}

if $q\sfo x\in a$ then $q\fo x\in b\,;$ \hfill 
if $q\sfo y\in b$ then $q\fo y\in a\,\msur$.

\itla{f-in}\msur
$p\fo b\in a$ \ iff \ 
$\forall\,q\<p\;\exists\,r\<q\;\exists\,z\;
(r\sfo z\in a \,\hbox{ and }\,r\fo b=z)$.
\een
%
Items \ref{f1} through \ref{f-in} define the forcing 
for formulas $a=b$ and $a\in b$ by induction on the ranks 
$\nrk a$ and $\nrk b$ of ``names'' $a,\,b\in\cN.$ The following 
items handle the standardness predicate and non--atomic formulas.
\ben
\def\theenumi{F-\arabic{enumi}}
\addtocounter{enumi}{3}
\itla{f-st}
$p\fo\st a$ \ iff \ 
$\forall\,q\<p\;\exists\,r\<q\;\est s\;(r\fo a=\brs)$.

\itla{f-neg}
$p\fo\neg\;\Phi$ \ iff \ none of $q\<p$ forces $\Phi$.

\itla{f-and}
$p\fo (\Phi\cj \Psi)$ \ iff \ $p\fo\Phi$ and $p\fo\Psi$.

\itla{f-all}
$p\fo \forall\,x\:\Phi(x)$ \ iff \ $p\fo\Phi(a)$ for every 
$a\in\cN$.
\een
It is assumed that the other logic connectives are 
combinations of ${\neg},\,\mathord{\cj},\,\forall.$

\ble
\label{monot}
Let\/ $a,\,b$ be ``names'' in $\cN.$  
If $p\sfo b\in a$ then ${p\fo b\in a}.$
 
If\/ $p\sfo b\in a$ and\/ $q\<p$ then\/ $q\sfo b\in a$.
 
If\/ $p\fo \Phi$ and\/ $q\<p$ then\/ $q\fo \Phi$.
\ele
\proof The first two assertions are quite obvious, the third one 
can be  easily proved by induction on the complexity of 
$\Phi$.\qed

\ble
\label{ok}
If\/ $p\in\dP$ does not force\/ $\Phi,$ a closed\/ 
\ste-formula with ``names'' in\/ $\cN$ as parameters, then 
there exists\/ $q\<p$ such that $q\fo\neg\;\Phi$.
\ele
\proof 
%
Assume ${\neg\;p\fo b\in a}.$ There exists a condition 
$q\<p$ such that 
$\neg\;\exists\,r\<q\;\exists\,z\;(r\sfo z\in a \cj r\fo b=z).$ 
To see that $q\fo\neg\;b\in a,$ let, on the contrary, a condition 
$q'\<q$ satisfy $q'\fo b\in a.$ Then by definition we have 
$r\sfo z\in a$ and $r\fo b=z$ for a condition $r\<q'$ and a 
``name'' $z,$ a contradiction with the choice of $q$.

Assume that ${\neg\;p\fo a=b}.$ Then by definition there exists  
$q'\<p$ such that e. g. for a ``name'' $x,$ $q'\sfo x\in a$ but 
$\neg\;q'\fo x\in b.$ 
It follows, by the above, 
that a condition $q\<q'$ satisfies $q\fo\neg\;x\in b.$ 
We prove that $q\fo a\not=b.$ Suppose that on the contrary a 
condition $r\<q$ forces $a=b.$ Since $r\sfo x\in a$ by 
Lemma~\ref{monot}, we have $r\fo x\in b,$ 
contradiction. 

A similar reasoning proves the result for formulas $\st a$.


As for non--atomic formulas, the result can be achieved by a 
simple straightforward induction on the logical complexity of 
the formula.
\qed

\subsection{Truth lemma}
\label{trul}

Suppose that $\Phi$ is a \ste-formula having ``names'' in $\cN$ 
as parameters. We let $\Phi[G]$ denote the formula obtained by 
\index{zzphig@$\Phi[G]$}
replacing occurrences of $\in$ and $\st$ in $\Phi$ by $\ing$ 
and $\stg,$ and every 
``name'' $a\in\cN$ by $a[G]\,;$  
thus $\Phi$ is a formula having sets in $\dH[G]$ as parameters.

\bte
\label{truth}
{\rm (The truth lemma.)} \ Let \/ $G\sq\dP$ be a\/ 
\dd\dP generic set over\/ $\dH.$
Let\/ $\Phi$ be a\/ \ste-formula having ``names'' in\/ $\cN$ as 
parameters. Then\/ $\Phi[G]$ is true in\/ $\dH[G]$ iff\/ 
$\exists\,p\in G\,(p\fo \Phi)$.
\ete
\proof Let us prove the result for the atomic formulas $a=b$ and 
$b\in a$ by induction on the ranks $\nrk$ of $a$ and $b.$ First 
of all we summarize the definition of the membership in $\dH[G]$ 
as follows: for all $a,\,b\in\cN$,
\dm
b[G]\ing a[G] \hspace{5mm}\hbox{iff}\hspace{5mm}
\exists\,b'\in\cN\;\exists\,p\in G\;
(b'[G]=b[G]\cj p\sfo b'\in a)\,.\eqno{(\ast)}
\dm 
%

We verify that $a[G]=b[G]$ iff some $p\in G$ satisfies $p\fo a=b.$  
Suppose that none of $p\in G$ forces $a=b.$ By the 
genericity, $G$ contains a condition $q$ such that, 
say, $q\sfo x\in a$ but $q\fo x\not\in b$ for some $x\in\cN.$ 
Then $x[G]\ing a[G]$ but 
$x[G]\ning b[G]$ by the induction hypothesis.  


Suppose now that $a[G]\not=b[G].$ Then, since $\dH[G]$ satisfies 
Extensionality, the sets differ from each other in $\dH[G]$ by 
their elements, say $x[G]\ing a[G]$ but $x[G]\ning b[G]$ for a 
``name'' $x\in\cN.$ By the induction hypothesis and $(\ast)$ there 
exist: a condition $p\in G$ and a ``name'' $x'$ such that 
$p\sfo x'\in a$ but $p\fo x'\not\in b.$ Then 
$p\fo a\not=b,$ because otherwise there exists a condition 
$q\<p$ which forces $a=b,$ immediately giving a contradiction. 

Consider a formula of the form $b\in a.$ Let a condition $p\in G$ 
force $b\in a.$ Then by the genericity of $G$ there 
exists a condition $r\in G$ such that 
$r\sfo z\in a$ and $r\fo z=b$ for a ``name'' $z\in\cN.$ 
This implies ${z[G]\ing a[G]}$ by definition and  
$z[G]=b[G]$ by the induction hypothesis. 

Assume now that ${b[G]\ing a[G]}$ and prove that a condition 
$p\in G$ forces $b\in a.$ We observe that, by $(\ast),$ there 
exist: a condition $p\in G$ and a ``name'' $b'$ such that 
$b'[G]=b[G]$ and $p\sfo b'\in a.$ We can assume that 
$p\fo b=b',$ by the induction hypothesis. Then  
$p\fo b\in a$ by definition. 

Formulas of the form $\st a$ are considered similarly. 
Let us proceed with non--atomic formulas by 
induction on the complexity of the formula involved. 

{\it Negation\/}. Suppose that $\Phi$ is $\neg\;\Psi.$ If 
$\Phi[G]$ is true then $\Psi[G]$ is false in $\dH[G].$ 
Thut none of $p\in G$ can force $\Psi,$ by the induction 
hypothesis. However the set 
$\ans{p\in\dP:p\,\hbox{ decides }\,\Psi}$ is dense in $\dP$ and 
belongs to $\dH.$ Thus some $p\in G$ forces 
$\Phi$ by the genericity of $G$. 

If $p\in G$ forces $\Phi$ then none of $q\in G$ can 
force $\Psi$ because $G$ is pairwise compatible. Thus  
$\Psi[G]$ fails in $\dH[G]$ by the induction hypothesis. 

{\it Conjunction\/}. Very easy. 

{\it The universal quantifier\/}. Let $p\in G$ force  
$\forall\,x\,\Psi(x).$ By definition we have $p\fo\Psi(a)$ for 
all $a\in\cN.$ Then $\Psi(a)[G]$ holds in $\dH[G]$ by the 
induction hypothesis, for all $a\in \cN.$ 
However $\Psi(a)[G]$ is $\Psi[G](a[G]),$ and 
$\dH[G]=\ans{a[G]:a\in\cN}.$ 
It follows that $\dH[G]\models\forall\,x\,\Psi(x)[G]$.

Assume that $\forall\,x\,\Psi(x)[G]$ is true in 
$\dH[G].$ By Lemma~\ref{ok} and the genericity some $p\in G$ 
forces either $\forall\,x\,\Psi(x)$ or $\neg\;\Psi(a)$ for a 
particular $a\in\cN.$ In the ``or'' case 
$\Psi(a)[G]$ is false in $\dH[G]$ by the induction hypothesis, 
contradiction with the assumption. Thus $p\fo\forall\,x\,\Psi(x),$ 
as required.\qed

\subsection{The extension models $\protect\HST$}
\label{hg:hst}

We complete the proof of Theorem~\ref{t:hg} in this subsection. 
Since the standard (therefore also internal) sets in $\dH[G]$ 
were already proved to be the same as in $\dH,$ we have the axioms 
of group~\ref{it1} (see Subsection~\ref{hst:ax}) in $\dH[G]$. 

Let us verify the $\ZFC$ axioms of group~\ref{hzfc} 
in $\dH[G].$ We concenrate on the axioms  
of Separation and Collection; the rest of the axioms can be easily 
proved following the $\ZFC$ forcing patterns. (Extensionality has 
already been proved, see Proposition~\ref{is}.) 

{\it Separation\/}. Let $X\in\cN,$ and $\Phi(x)$ be a \ste-formula 
which may contain sets in $\cN$ as parameters. We have to find a 
``name'' $Y\in \cN$ satisfying $Y[G]=\ans{x\in X[G]:\Phi[G](x)}$ 
in $\dH[G].$ Note that by definition all 
elements 
of $X[G]$ in $\dH[G]$ are of the form $x[G]$ where $x$ belongs to 
the set $\cX=\ans{x\in\cN:\exists\,p\;(\ang{p,x}\in X)}\in\dH.$ 
(We suppose that $\nrk X>0.$ The case $X\in\cN_0$ does not differ 
much.)
Now $Y=\ans{\ang{p,x}\in \dP\ti \cX:p\fo \Phi(x)}$ is the required 
``name''. (See Shoenfield~\cite{sh} for details.)

{\it Collection\/}. Let $X\in\cN,$ and $\Phi(x,y)$ be a formula  
with ``names'' in $\cN$ as parameters. We have to find a ``name'' 
$Y\in\cN$ such that
\dm
\forall\,x\in X[G]\;(\exists\,y\;\Phi[G](x,y)\;\lra\;
\exists\,y\in Y[G]\;\Phi[G](x,y))
\dm
is true in $\dH[G].$ Let $\cX\in\dH,\;\cX\sq\cN$ be defined as above,  
in the proof of Separation. Using Collection in $\dH,$ we obtain a 
set $\cY\sq\cN,$ 
sufficient in the following sense: if $x\in \cX,$ and $p\in \dP$ 
forces $\exists\,y\;\Phi(x,y),$ then 
\dm
\forall\,q\<p\;\exists\,r\<q\;\exists\,y\in \cY\;
(r\fo\Phi(x,y))\,.
\dm
The set $Y=\dP\ti \cY$ (then $Y[G]=\cY$) is as required.

{\it Weak Regularity\/}. Let $X\in\cN.$ Using an appropriate 
dense set in $\dP,$ we find a condition $p\in G$ such that 
$p\fo a\in X$ for a ``name'' $a\in\cN,$ but 
{$1)\msur$ $p\fo b\not\in X$} for any ``name'' $b\in\cN_\ba$ 
where $\ba<\nrk a$ -- provided $\nrk a>0,$ and 
{$2)\msur$ $p\fo \bry\not\in X$} for any $y\in\dH$ with 
$\irk y<\irk x$ -- provided $a=\brx\in\cN_0.$ Now, if 
$\nrk a>0,$ or if $a=\brx\cj \irk x>0,$ then simply 
$p\fo a\cap X=\emps.$ If finally $a=\brx$ for an internal $x$ 
then $x\cap X[G]$ contains only internal elements.

{\it Standardization\/}. Let $X\in\cN.$ We have to find a standard 
set $Y$ which contains in $\dH[G]$ the same standard elements as 
$X[G]$ does. It can be easily proved by induction on $\nrk a$ that, 
for every ``name'' $a\in\cN$,   
\dm
\stan a=\ans{s:\st s\cj\exists\,p\in\dP\;(p\fo\brs\in a)}
\dm
is a set in $\dH.$ 
($\stan a$ contains all standard \dd\ing elements of $a[G]$.) 
Thus $\stan X\sq S$ for a standard $S,$ by Lemma~\ref{hstb}. 
Since $\dP$ is standard size distributive in $\dH,$ $G$ contains, by 
the genericity, a condition $p$ which, for any standard $s\in S,$ 
decides the statement $\brs\in X.$ Applying Standardization in 
$\dH,$ we get a standard set $Y\sq S$ such that, for each 
$s\in S,$ $s\in Y$ iff $p\fo\brs\in X.$ The $Y$ is as required. 

{\it Standard size Choice\/}. The problem can 
be reduced to the following form. Let $S$ be a standard set,  
$P\in\cN,$ and $P[G]$ is a set of pairs in $\dH[G].$  
Find a ``name'' $F$ such that the following is true in $\dH[G]:$ 
\its 
\bit
\item[]\msur
{\it $F[G]$ is a function defined on\/ $\upsG S$ and satisfying\/\\  
${\exists\,y\,P[G](x,y)\;\lra\;P[G](x,F[G](x))}$ for each standard\/ 
$x\in S$.}\its
\eit
Arguing as above and using the standard size Choice in 
$\dH,$ we obtain a condition $p\in G$ and a 
function $f\in\dH,$ $f:\upsG S\;\lra\;\cN,$ such that, for every 
$x\in\upsG S,$ $p$ either forces $\neg\;\exists\,y\;P(\brx,y)$ or 
forces $P(\brx,y_x)$ where $y_x=f(x)\in\cN.$ 
One easily converts $f$ to a required ``name'' $F$.

{\it Dependent Choice\/} -- similar reduction to $\dH$.

{\it Saturation\/}. Using the same argument, one proves that each  
standard size family of internal sets in $\dH[G]$ already belongs 
to $\dH.$ \qed

\np

%
%
\section{Generic isomorphisms}
\label{isom}

Let us consider a particular forcing which leads to a 
generic isomorphism between two internally presented 
elementarily equivalent structures. 

We continue to consider a model $\dH\models\HST$ satisfying 
assumptions \ref{dhis} and \ref{dhwf} in Section~\ref{f}. 
We suppose in addition that\its

\ben
\def\theenumi{$(\fnsymbol{enumi})$}
\def\labelenumi{\theenumi}
\addtocounter{enumi}{3}
\itla{lab}\msur 
$\cL\in\dH$ is a first--order language containing 
(standard size)--many symbols. $\gA=\ang{A;...}$ and 
$\gB=\ang{B;...}$ are two internally presented elementarily 
equivalent \dd\cL structures in $\dH$.\its
\index{structures!ab@$\gA$ and $\gB$}
\index{language!l@$\cL$}
\een
By definition both $A$ and $B$ are internal sets, and the 
interpretations of 
each symbol of $\cL$ in $\gA$ and $\gB$ are internal in $\dH$. 

The final aim is to obtain a generic isomorphism $\gA$ onto $\gB.$ 
We shall define a notion of forcing $\dP=\dP_{\cL\gA\gB}\in\dH$  
such that $\gB$ is isomorphic to $\gA$ in every \dd\dP generic 
extension of $\dH,$ provided $\dH$ satisfies a requirement which 
guarantees the standard size distributivity of $\dP$. 

It is the most natural idea to choose the forcing conditions 
among partial 
functions $p,$ mapping subsets of $A$ onto subsets of $B.$ We 
have to be careful: the notion of forcing must be a set, thus 
for instance maps having a standard size domain do not work 
because even an \dd\dI finite internal infinite set has a proper 
class of standard size subsets in $\HST.$ We overcome this 
obstacle using {\it internal\/} partial maps $p,$ such 
that each $a\in\dom p$ satisfies in $\gA$ exactly the same 
\dd\cL formulas as $p(a)$ does in $\gB$. 

Given a condition $p$ and an element 
$\fa\in A\setminus\dom p,$ we must be able to incorporate $\fa$ 
in $p,$ i.\ e. define a stronger 
condition $p_+$ such that $\fa\in\dom p_+.$ Here we face a 
problem: to find an element $\fb\in B$ which, for each $a\in\dom p,$ 
is in the same relations with $p(a)$ in $\gB$ as $\fa$ is with 
$a$ in $\gA.$ Since $\dom p$ cannot be a set of standard size, 
it is not immediately  clear how a saturation argument can be used 
to get a required $\fb$. 

We shall develop the idea as follows. Let $\Phi(x,y)$ be an 
\dd\cL formula. We are willing to find $\fb\in B$ so that 
$\Phi(\fa,a)$ in $\gA$ iff $\Phi(\fb,p(a))$ in $\gB$ for all 
$a\in\dom p.$ The sets $u=\ans{a\in\dom p:\gA\mo\Phi(\fa,a)}$ 
and $v=\dom p\setminus u$ are internal by the choice of $\gA.$ 
We observe that the chosen element $\fa$ satisfies 
$\forall\,a\in u\;\Phi(\fa,a)$ and 
$\forall\,a\in v\;\neg\,\Phi(\fa,a)$ in $\gA,$ so that the 
sentence 
\dm
\exists\,x\;[\,\forall\,a\in u\;\Phi(x,a)\cj 
\forall\,a\in v\,\neg\,\Phi(x,a)\,]
\dm 
is true in $\gA.$ Suppose that $p$ 
also preserves sentences of this form, so  
that 
\dm
\exists\,y\;[\,\forall\,b\in p\ima u\;\Phi(y,b)\cj 
\forall\,b\in p\ima v\,\neg\,\Phi(y,b)\,]
\dm
is true in $\gB.$ ($p\ima u=\ans{p(a):a\in u}$ is the \dd pimage 
of $u$.) This gives an element $\fb\in B$ which may be 
put in correspondence with $\fa$.

Thus we have to preserve formulas of the displayed type, i.\ e. 
\dd\cL formulas with some internal subsets of $A$ as parameters, 
so in fact a stronger preservation hypothesis is involved than 
the result achieved. This leads to a sort of hierarchical 
extension of the language $\cL$. 

\subsection{The extended language}

Arguing in $\dH,$ we define the notion of {\it type\/} as follows. 
Let $D\in\dI$. 

$0$ is a type. An object of type $0$ over a set $D$ 
is an element of $D$. 

Suppose that $l_1,...,l_k$ are types. Then $l=\tp(l_1,...,l_k)$ is 
a type. (Here $\tp$ is a formal sign.) An object of type $l$ over 
$D$ is an internal set of \dd ktuples $\ang{x_1,...,x_k}$ where 
each $x_i$ is an object of type $l_i$ over $D$.  
\index{type}
\index{zztau@$\tp(l_1,...,l_k)$}
\index{zzli@$\cLi$}

E. g. objects of type $\tp(0,0)$ over $D$ are internal 
subsets of $D\ti D$.

We define $\cLi$ as the extension of $\cL$ by variables 
\index{language!le@extended, $\cLi$}
$x^l,y^l,...$ for each type $l,$ which can enter 
formulas~\footnote
{\ By {\it formulas\/}, with respect to the languages $\cL$ and 
$\cLi,$ we shall understand finite sequences of something 
satisfying certain known requirements, not only 
``metamathematical'' formulas. Since any finite tuple of internal 
sets is internal (an easy consequence of Lemma~\ref{exten}), a 
formula with internal parameters is formally an internal object, 
so for instance its truth domain is internal as well.} 
only through the expressions $x^l(x^{l_1},...,x^{l_k})$ (may 
be written as $\ang{x^{l_1},...,x^{l_k}}\in x^l$), provided 
$l=\tp(l_1,...,l_k),$ and also $x=x^0,$ where $x$ is an 
\dd\cL variable. (We shall formally distinguish variables of 
type $0$ from \dd\cL variables.) 

Let $\gC=\ang{C;...}$ be an internally presented \dd\cL structure. 
Given an internal set $D\sq C,$ we define a type--theoretic 
extended structure $\gC[D]$ which includes the ground domain $C$ 
\index{zzcd@$\gC[D]$}
with all the \dd\gC interpretations of \dd\cL symbols, and the domain 
$D^l=\ans{x^l:x^l\,\hbox{ is an object of type }\,l\,
\hbox{ over }\,D}$ 
\index{zzdl@$D^l$}
for each~type~$l.$ 

We observe that each $D^l$ is an internal set 
because the construction of $D^l$ can be executed in $\dI.$  
For instance $D^{\tp(0)}=\cP(D)$ in $\dI.$ 
$D=D^0$ is an internal subset of $C$. 

Every \dd\cLi formula (perhaps, containing sets in 
$\gC[D]$ as parameters) can be interpreted in $\gC[D]$ in the 
obvious way. (Variables of type $l$ are interpreted in $D^l.$) 
This converts $\gC[D]$ to an internally presented  
\dd\cLi structure.

It will always be supposed that $D^{l_1}\cap D^{l_2}=\emps$ 
provided $l_1\not=l_2$. 

We put $\Di=\bigcup_l D^l$.
\index{zzdi@$\Di$}

\subsection{The forcing}
\label{if}

We recall that a standard size language $\cL$ and a pair of 
internally presented elementarily equivalent \dd\cL models 
$\gA=\ang{A;...}$ and $\gB=\ang{B;...}$ are fixed. 

Suppose that $p$ is an internal $1-1$ map from an internal set 
$D\sq A$ onto a set $E\sq B.$ We expand $p$ on all types 
$l$ by induction, putting 
\dm
p^l(x^l)=\ans{\ang{p^{l_1}(x^{l_1}),...,p^{l_k}(x^{l_k})}:
\ang{x^{l_1},...,x^{l_k}}\in x^l}
\hspace{3mm}\hbox{for all}\hspace{3mm}x^l\in D^l,
\dm 
whenever $l=\tp(l_1,...,l_k).$ Then $p^l$ 
\index{zzpl@$p^l(x^l)$}
internally $1-1$ maps $D^l$ onto $E^l$. 

If $\Phi$ is an \dd\cLi formula containing parameters in $\Di$  
then let $p\Phi$ be the formula obtained by changing each 
\index{zzpphi@$p\Phi$}
parameter $x\in D^l$ in $\Phi$ to $p^l(x)\in E^l.$

\bdf
\label{iforc}
$\dP=\dP_{\cL\gA\gB}$ is the set of all internal $1-1$ maps $p$ 
\index{zzdplab@$\dP_{\cL\gA\gB}$}
such that $D=\dom p$ is an (internal) 
subset of $A,$ $E=\ran p\sq B$ (also internal), and, for each 
closed \dd\cLi formula $\Phi$ having sets in $\Di$ as 
parameters, we have $\gA[D]\mo\Phi$ iff $\gB[E]\mo p\Phi$. 

We define $p\<q$ ($p$ is stronger than $q$) iff $q\sq p$.\qed
\edf
For instance the empty map $\emptyset$ belongs to $\dP$ because 
$\gA$ and $\gB$ are elementarily equivalent. (The properly 
\dd\cLi variables can be eliminated in this case because the 
domains become finite.) 
We shall see 
%
(Corollary~\ref{t:isom}, the main result) that the forcing leads to 
generic isomorphisms $\gA$ onto $\gB.$ This is based on the following 
two technical properties of this forcing.

\bpro
\label{fo1}
$\dP=\dP_{\cL\gA\gB}$ is standard size closed in $\dH$.
\epro 

\bpro
\label{fo2}
Let\/ $p\in\dP,\;D=\dom p,\;E=\ran p.$ If\/ 
$\fa\in A\setminus D$ then there exists\/ $\fb\in B\setminus E$ 
such that $p_+=p\cup\ans{\ang{\fa,\fb}}\in\dP.$ Conversely, if\/ 
$\fb\in B\setminus E$ then there exists\/ $\fa\in A\setminus D$ 
such that $p_+=p\cup\ans{\ang{\fa,\fb}}\in\dP$.
\epro

\proof{}of Proposition~\ref{fo1}. Let $\kpa$ be a cardinal. 
Suppose that $p_\al\;\,(\al<\kpa)$ are conditions in $\dP,$ and 
$p_\ba\<p_\al$ whenever $\al<\ba<\kpa.$ By definition $\dP$ is a 
standard size intersection of internal sets (because the structures 
are internally presented and $\cLi$ is a language of standard size), 
hence so is each of the sets $P_\al=\ans{p\in\dP:p\<p_\al}.$ 
Furthermore $P_\al\not=\emps$ and 
$P_\ba\sq P_\al$ whenever $\al<\ba<\kpa.$ Finally $\kpa$ 
(as every set in the condensed universe~$\dV$) is a set of standard 
size by Lemma~\ref{wo=ss}, so $\bigcap_{\al<\kpa}P_\al\not=\emps$ 
by Saturation. Thus there exists $p\in P$ such that $p\<p_\al$ for 
all $\al,$ as required.\qed

\bdf
\label{loc}
A function $F$ defined on an internal set $A$ is 
{\it locally internal\/} iff either $A$ is finite or for any 
\index{locally internal map}
$a\in A$ there exists an infinite internal set $A'\sq A$ 
containing $a$ and such that $F\res A'$ is internal.\qed
\edf

\bcor
\label{t:isom}
Suppose that, in addition to\/ \ref{dhis}, \ref{dhwf}, \ref{lab}, 
$\dH$ satisfies statement~\ref{gdc} in Theorem~\ref{li:hst}. 
Let\/ $\dP=\dP_{\cL\gA\gB}$ in\/ $\dH.$ Then, every\/ 
\dd\dP generic extension\/ $\dH[G]$ is a model of\/ $\HST,$ a 
plain extension of\/ $\dH,$ and $F=\bigcup G$ is a locally 
internal isomorphism\/ $\gA$ onto\/ $\gB$ in $\dH[G]$.
\ecor
\proof{}of the corollary. Proposition~\ref{fo1} plus the assumed 
statement~\ref{gdc} guarantee that $\dP$ is standard size 
distributive in $\dH$ in the sense of Definition~\ref{ssdis}. 
Therefore $\dH[G]\mo\HST,$ by Theorem~\ref{t:hg}. Furthermore 
$F$ $1-1$ maps $A$ onto $B$ by Proposition~\ref{fo2}. The map is 
an isomorphism because $F$ is a union of conditions $p\in\dP$ which 
preserve the truth of \dd\cL sentences. 

To prove that $F$ is locally internal, assume that the sets $A,\,B$ 
are infinite, and $a\in A.$ It can be easily proved by the same 
reasoning as in the proof of Proposition~\ref{fo1} that $G$ 
contains a condition $p$ such that $a\in\dom p$ and $\dom p$ is 
infinite (although perhaps \dd\dI finite).\qed\eproof

The remainder of this section is devoted to the proof of 
Proposition~\ref{fo2}. By the simmetry, we concentrate on the first 
part. Let us fix a condition $p\in\dP.$ Let $D=\dom p,\;E=\ran p.$ 
\index{zzapde@$\fa,\,p,\,D,\,E$}
(For instance we may have $p=D=E=\emps$ at the moment.)
Consider an arbitrary $\fa\in A\setminus D;$ we have to find a 
counterpart  $\fb\in B\setminus E$ such that 
$p_+=p\cup\ans{\ang{\fa,\fb}}\in\dP$. 

\subsection{Adding an element}
\label{adding}

Let $\kpa=\card\cL$ (or $\kpa=\aleph_0$ provided $\cL$ is finite) 
in $\dH.$  
We enumerate by $\Phi_\al(x)\;(\al<\kpa)$ all parameter--free 
\index{zzphial@$\Phi_\al$}
\dd\cLi formulas which contain only one \dd\cL variable $x$ but 
may contain several variables $x^l$ for various types $l$.

Let us consider a particular \dd\cLi formula 
$\Phi_\al(x)=\vpi(x,x^{l_1},...,x^{l_n}).$ Let $l=\tp(l_1,...,l_n).$ 
(Both $l$ and each of $l_i$ are types.) We set 
\dm
X_\al=  
\ans{\ang{x^{l_1},...,x^{l_n}}\in D^{l_1}\ti...\ti D^{l_n}:
\gA[D] \mo\vpi(\fa,x^{l_1},...,x^{l_n})}\,;
\dm
\index{zzxal@$X_\al$}
thus $X_\al$ is internal~\footnote
{\ To prove that $X_\al$ is internal, we first note that every 
finite subset of $\dI$ is internal, which can be easily proved 
by induction on the number of elements. Therefore, since only 
finitely many types $l$ are actually involved in the definition 
of $X_\al,$ while all the relevant domains and relations are 
internal, the definition of $X_\al$ can be executed in $\dI.$ 
(We cannot, of course, appeal to Transfer since $\vpi$ is not a 
metamathematical formula here.)}
and $X_\al\in D^l.$ Let $\Psi_\al(X_\al,x)$ be the \dd\cLi formula 
\index{zzpsial@$\Psi_\al$}
\dm
\forall\,x^{l_1}\,...\,\forall\,x^{l_n}\;[\,
X_\al(x^{l_1},...,x^{l_n})\; \llra \;
\vpi(x,x^{l_1},...,x^{l_n})\,]\,,
\dm 
so that by definition $\gA[D]\models \Psi_\al(X_\al,\fa).$ 
Thus we have \dd\kpa many formulas $\Psi_\al(X_\al,x),$ realized 
in $\gA[D]$ by one and the same element $x=\fa\in A$. 

We put $Y_\al=p^l(X_\al);$ so that $Y_\al\in E^l$.
\index{zzyal@$Y_\al$}

\ble
\label{eb}
There exists\/ $\fb\in B$ which realizes (by\/ $y=\fb$) every 
formula\/ $\Psi_\al(Y_\al,y)$ in $\gB[E]$.
\ele
\proof By Saturation, it suffices to prove that every finite 
conjunction 
$\Psi_{\al_1}(Y_{\al_1},y)\cj ... \cj \Psi_{\al_m}(Y_{\al_m},y)$ 
can be realized in $\gB[E].$ By definition $\fa$ witnesses that 
${\gA[D]\mo\exists\,x\;[\,\Psi_{\al_1}(X_{\al_1},x)\cj ... 
\cj \Psi_{\al_m}(X_{\al_m},x)\,]\,}.$ Therefore  
$\gB[E]\mo\exists\,y\;[\,\Psi_{\al_1}(Y_{\al_1},y)\cj ... 
\cj \Psi_{\al_m}(Y_{\al_m},y)\,]\,,$ since $p\in\dP$.\qed\eproof

Let us fix an element  $\fb\in B$ satisfying $\Psi_\al(Y_\al,\fb)$ 
in $\gB[E]$ for all $\al<\kpa.$ We set 
$p_+=p\cup\ans{\ang{\fa,\fb}},\;\,D_+=D\cup\ans{\fa},\;\,
E_+=E\cup\ans{\fb}$. 
\index{zzbpde@$\fb,\,p_+,\,D_+,\,E_+$}

\subsection{Why the choice is correct}

It will take some effort to check that $p_+$ is a condition in 
$\dP.$ Let 
us prove first a particular lemma which shows that $p_+$ preserves 
formulas containing $\fa$ and sets in $\Di$ as parameters. 

\ble
\label{+1}
Let\/ $\vpi(x)$ be an\/ \dd\cLi formula which may contain sets 
in\/ $\Di$ as parameters. Then\/ $\vpi(\fa)$ is true in\/ $\gA[D]$ 
iff\/ $(p\vpi)(\fb)$ is true in\/ $\gB[E]$.
\ele
\proof The formula $\vpi(x)$ is obtained from a parameter--free 
\dd\cLi formula $\Phi(x,\dots)$ by changing free variables in the 
list $\dots$ to appropriate parameters (of the same type) from 
$\Di.$ We can assume that in fact the list $\dots$ does not 
include \dd\cL variables; indeed if such one, say $y,$ occurs 
then we first change $\Phi(x,y,...)$ to 
$\exists\,y\;[\,\Phi(x,y,...)\cj y=y^0\,]\,,$ where $y^0,$ a 
variable of type $0,$ is free. In this assumption, 
$\vpi(x)$ is $\Phi_\al(x,x^{l_1},...,x^{l_n})$ for some $\al$ and 
parameters $x^{l_i}\in D^{l_i}.$ 
Then, since $\Psi_\al(X_\al,\fa)$ is true in $\gA[D],$ we have
\dm
X_\al(x^{l_1},...,x^{l_n})
\hspace{5mm}\hbox{iff}\hspace{5mm}
\gA[D]\mo\Phi_\al(\fa,x^{l_1},...,x^{l_n}).
\dm
Note that 
$X_\al(x^{l_1},...,x^{l_n})\,\llra\,Y_\al(y^{l_1},...,y^{l_n}),$ 
where $y^{l_i}=p^{l_i}(x^{l_i})\in E^{l_i},$ since $p\in\dP.$ 
Finally,  $Y_\al(y^{l_1},...,y^{l_n})$ iff 
$\gB[E]\mo\Phi_\al(\fb,y^{l_1},...,y^{l_n}),$ 
because $\Psi_\al(Y_\al,\fb)$ is true in $\gB[E]$ by the choice of 
$\fb.$ However the formula $\Phi_\al(\fb,y^{l_1},...,y^{l_n})$ 
coincides with $(p\vpi)(\fb)$.\qed\eproof

Taking the formula $x\not\in D$ as $\vpi(x)$ (then $p\vpi(x)$ is 
$x\not\in E$), we obtain $\fb\not\in E,$ so $p_+$ is a 
$1-1$ internal map. It remains to check that $p_+$ transforms true 
\dd\cLi formulas with parameters in $\Dpi$ into true 
\dd\cLi formulas with parameters in $\Epi$. The idea is 
to convert a given formula with parameters in $\Dpi$ into 
a \dd\cLi formula with parameters in $\Di$ plus $\fa$ as an extra 
parameter, and use Lemma~\ref{+1}. 

Fortunately the structure of types over an internal set $C$ depends 
only on the internal cardinality of $C$ but does not depend on the 
place $C$ takes within $\gA.$ This allows to ``model'' $\Dpi$ 
in $\Di$ identifying the $\fa$ with $\emps$ and any $a\in D$ 
with $\ans{a}.$ To realize this plan, let us define 
$
\cD=\ans{\emps}\cup\ans{\ans{a}:a\in D}\,,
\index{zzdc@$\cD$}
$
so that $\cD\sq D^\ell,$ where $\ell=\tp(0)$ (the type of 
\index{zzlell@$\ell$}
subsets of $D$). Furthermore we have $\cD\in D^{\tp(\ell)}$ 
because $\cD$ is internal. 

For each type $l,$ we define a type $\da(l)$ by $\da(0)=\ell$ 
\index{zzdl@$\da(l)$}
and $\da(l)=\tp(\da(l_1),...,\da(l_n))$ provided 
$l=\tp(l_1,...,l_n)$.

We put $\da(\fa)=\emps,$ and $\da(a)=\ans{a}$ for 
all $a\in D,$ so that $\da$ is an internal bijection $D_+$ 
onto $\cD.$ We expand $\da$ on higher types by 
$\da(x)=\ans{\ang{\da(x_1),...,\da(x_n)}:\ang{x_1,...,x_n}\in x};$ 
\index{zzdx@$\da(x)$}
thus $\da(x)\in \cD^l\sq D^{\da(l)}$ whenever $x\in {D_+}^l.$ 
Take notice that $\da({D_+}^l)=\cD^l.$ Thus 
$\da=\da_{D\fa}$ defines a $1-1$ correspondence between 
\index{zzdad@$\da_{D\fa}$}
$\Dpi$ and $\Dsi$. 

Let $\psi(x^{l_1},...,x^{l_n},v_1,...,v_m)$ be a 
\dd\cLi formula, containing \dd\cL variables $v_j$ and 
properly \dd\cLi variables $x^{l_i}.$ We introduce 
another \dd\cLi formula, denoted by  
$\psi_{D\fa}
(\xi^{\da(l_1)},...,\xi^{\da(l_n)},v_1,...,v_m),$ 
\index{zzpsid@$\psi_{D\fa}$}
containing $\fa,\,D,$ and 
finitely many sets $\cD^l$ as parameters (this is symbolized 
by the subscript $D;$ the involved sets $\cD^l$ are derivates 
of $D$ and $\fa$), as follows. 

Each free variable $x^{l_i}$ is changed to some $\xi^{\da(l_i)},$ 
a variable of type $\da(l_i).$ (We use characters $\xi,\,\eta,\,\za$ 
for variables intended to be restricted to $\Dsi$). 

Each quantifier $\rQ\,u^l\;...\,u^l\,...$ 
is changed to 
$\rQ\,\eta^{\da(l)}\in \cD^l\,...\,\eta^{\da(l)}...\;.$ 
(Note that 
$\cD^l=\da({D_+}^l)$ is an internal subset of $D^{\da(l)}.$)

Each occurrence of type $x=\za^\ell$ (which is obtained by the 
abovementioned transformations from an original equality $x=z^0$) 
is changed to 
\dm
(x=\fa\cj \za^\ell=\emps) \orr 
(x\in D\cj \za^\ell=\ans{x}) \eqno{(\ast)}
\dm
(the equalities $\za^\ell =...$ can here be converted to correct 
\dd\cLi formulas). 

\ble
\label'
Let\/ $\psi(x^{l_1},...,x^{l_n},v_1,...,v_m)$ be an\/ 
\dd\cLi formula, $x^{l_i}\in {D_+}^{l_i}$ and\/ $a_j\in A$ 
for all\/ $i$ and\/ $j.$ Then\/ 
$\gA[D_+]\mo\psi(x^{l_1},...,x^{l_n},a_1,...,a_m)$ iff\/ 
$\gA[D]\mo
\psi_{D\fa}(\da(x^{l_1}),...,\da(x^{l_n}),a_1,...,a_m)$.
\ele
\proof Suppose that $\psi$ is an atomic formula. If $\psi$ is 
in fact an \dd\cL formula then the equivalence is obvious. 
Otherwise $\psi$ is either of the form $x=x^0$ or of the form 
$\ang{x^{l_1},...,x^{l_n}}\in x^l,$ where $l=\tp(l_1,...,l_n).$ 
The latter case does not cause a problem: use the definition 
of $\da(x^l)$.

Consider a formula of the form $x=z^0$ as $\psi(z^0,x),$ where 
$x\in A$ and $z^0\in D_+={D_+}^0.$ By definition $\da(0)=\ell,$ 
$\da(z^0)=\emps$ provided $z^0=\fa$ and 
$\da(z^0)=\ans{z^0}$ otherwise, and 
$\psi_{D\fa}(\za^\ell,x)$ is the formula $(\ast).$ 
One easily sees that $x=z^0$ iff 
$\psi_{D\fa}(\da(z^0),x)$.


As for the induction step, we consider only the step $\exists\,u^l$ 
because the connectives $\neg$ and $\cj$ are automatical, as well as 
$\exists$ in the form $\exists\,x,$ where $x$ is an \dd\cL variable. 

Let $\psi(...)$ be the formula $\exists\,u^l\,\phi(u^l,...).$ We 
have the following chain
:
\ben
\def\theenumi{(\arabic{enumi})}
\def\labelenumi{\theenumi}

\itla{2)} $\gA[D_+]\mo\psi(...)$, \ that is, \ 
$\exists\,u^l\in {D_+}^l\,(\gA[D_+]\mo\phi(u^l,...))$

\itla{3)} $\exists\,u^l\in {D_+}^l\,(\gA[D]\mo
\phi_{D\fa}(\da(u^l),...))$

\itla{4)} $\exists\,\eta^{\da(l)}\in \cD^l\,(\gA[D]\mo
\phi_{D\fa}(\eta^{\da(l)},...))$

\itla{5)} $\gA[D]\mo\psi_{D\fa}(...)$
\een
The equivalence $\ref{2)}\,\llra\,\ref{3)}$ holds by induction 
hypothesis, equivalence $\ref{3)}\,\llra\,\ref{4)}$ follows from 
the equality $\cD^l=\ans{\da(u^l):u^l\in {D_+}^l}\sq D^{\da(l)},$ 
and the equivalence $\ref{4)}\,\llra\,\ref{5)}$ from the fact 
that $\psi_{D\fa}(...)$ is the formula 
$\exists\,\eta^{\da(l)}\in \cD^l\,
\phi_{D\fa}(\eta^{\da(l)},...)$ by definition.\qed

\subsubsection*{We complete the proof of 
Proposition~\protect\ref{fo2}} 

Let $\Phi$ be the \dd\cLi formula $\phi(x_1,...,x_n)$ 
containing sets $x_i\in {D_+}^{l_i}$ as parameters. Let 
$y_i=p^{l_i}(x_i);$ so that $y_i\in {E_+}^{l_i}.$ Let $\Psi$ be 
the \dd\cLi formula $\phi(y_1,...,y_n).$ We have to prove that \ 
$\gA[D_+]\mo\Phi$ \ iff \ $\gB[E_+]\mo\Psi$. \vspace{1mm}

{\it Step 1\/}. $\gA[D_+]\mo\Phi$ \ iff \ 
$\gA[D]\mo\phi_{D\fa}(\da(x_1),...,\da(x_n))$  
(Lemma~\ref').\vspace{1mm}

{\it Step 2\/}. Let  
$\cE=\ans{\emps}\cup\ans{\ans{b}:b\in E}=p\ima\cD.$ 
We observe that the final statement of 
step 1 is equivalent, by Lemma~\ref{+1}, to the following one: 
$\gB[E]\mo\phi_{E\fb}(p(\da(x_1)),...,p(\da(x_n)))$.%
\vspace{1mm} 

{\it Step 3\/}. 
In the last formula, $\da=\da_{D\fa}$ is 
the transform determined by $D$ and $\fa.$ Let us consider 
its counterpart, $\vep=\da_{E\fb}.$ One can easily verify 
that then $p(\da(x_i))=\vep(y_i),$ where, we recall, 
$y_i=p^l(x_i).$ So the final statement of step 2 is equivalent to 
$\gB[E]\mo\phi_{E\fb}(\vep(y_1),...,\vep(y_n))$.\vspace{1mm} 

{\it Step 4\/}. Using Lemma~\ref' with respect to the transform 
$\vep=\da_{E\fb}$ and the model $\gB,$ we conclude that the final 
statement of step 3 is equivalent to $\gB[E_+]\mo\Psi$.\qed 

\np

%
\section{A model for the isomorphism property} 
\label{total}

Fortunately the generic extensions of the considered type do 
not introduce new internal sets and new standard size collections  
of internal sets. This makes it possible to ``kill'' all pairs 
of elementarily equivalent internally presented structures by a 
product rather than iterated forcing. Following this idea, we 
prove (Theorem~\ref{conip} below) that a product $\Pi$ of 
different forcing notions of the form $\dP_{\cL\gA\gB},$ with 
internal \dd\dI finite support, leads to generic extensions which 
model $\HST$ plus the isomorphism property $\ip$. 

The product forcing will be a class forcing in this case because 
we have class--many pairs to work with; this will cause some 
technical problems in the course of the proof, in comparison 
with the exposition in Section~\ref f. 


We continue to consider a model $\dH$ of $\HST,$ satisfying 
requirements \ref{dhis} and \ref{dhwf} in Section~\ref{f}. 
$\dS\sq \dI$ and $\dV$ are resp. the classes of all 
standard and internal sets in $\dH,$ and the condensed subuniverse. 

\subsection{The product forcing notion}
\label{prod}

Arguing in $\dH,$ let us enumerate 
somehow all relevant triples consisting of a language $\cL$ and 
a pair of \dd\cL structures $\gA,\,\gB,$ to be made isomorphic. 

Let, in $\dH,$ $\Ind$ be the class of all \dd5tuples 
\index{zzind@$\Ind$}
$i=\ang{w,\kpa,\fL,\fA,\fB}$ such that $w$ is an internal  
set, $\kpa$ is an \dd\dI cardinal, $\fL=\ans{s_\al:\al<\kpa}$ a 
first--order internal language \dd\dI containing $\<\kpa$ symbols, 
and $\fA,\,\fB$ are internal \dd\fL structures. 
(Then obviously $i$ itself is internal.) 

We set $w_i=w,$ $\kpa_i=\kpa,$ $\fL_i=\fL,$ $\fA_i=\fA,$ $\fB_i=\fB$.
\index{zzwietc.@$w_i,\,\kpa_i,\,\fL_i,\,\fA_i,\,\fB_i$}

It is clear that $\Ind$ is a class \mem-definable in $\dI$. 
Elements of $\Ind$ will be called {\it indices\/}. 
\index{index}

Suppose that $i\in\Ind.$ Then by definition 
$\fL_i=\fL=\ans{s_\al:\al<\kpa}$ is an internal language (with a 
fixed internal enumeration of the \dd\fL symbols). 
We define the {\it restricted\/} standard size language  
$\cL=\cL_i=\ans{s_\al:\al<\kpa\cj\st\al}.$ 
\index{zzlietc@$\cL_i,\,\gA_i,\,\gB_i$}
Let $\gA_i$ and $\gB_i$ 
denote the corresponding restrictions of $\fA$ and $\fB;$ then 
both $\gA_i$ and $\gB_i$ are internally presented 
\dd{\cL_i}structures. 

On the other hand, if $\cL$ is a standard size language and 
$\gA,\,\gB$ a pair of internally presented \dd{\cL}structures 
then there exists an index $i\in\Ind$ such that $\cL=\cL_i,$ 
$\gA=\gA_i,$ and $\gB=\gB_i$.

The forcing $\Pi$ will be defined as a collection of internal 
functions $\pi.$ Before the exact definition is formulated, 
let us introduce a useful notation: $\avl\pi=\dom\pi$ (then 
\index{zzpiavl@$\avl\pi$}
$\avl\pi\in\dI$) and $\pi_i=\pi(i)$
\index{zzpii@$\pi_i$}
for all $\pi\in\Pi$ and $i\in\avl\pi$. 

\bdf
\label{totp}
$\Pi$ is the collection of all internal functions $\pi$ such 
\index{zzPi@$\Pi$}
\index{forcing!pi@$\Pi,$ product forcing}
that $\avl\pi\sq\Ind$ is an \dd\dI finite (internal) set, and 
$\pi_i\in\dP_{\cL_i\,\gA_i\,\gB_i}$ for each $i\in\avl\pi.$ We set 
$\pi\<\rho$ (i. e. $\pi$ is stronger than $\rho$) iff 
$\avl\rho\sq\avl\pi$ and $\pi_i\<\rho_i$ (in 
$\dP_{\cL_i\,\gA_i\,\gB_i},$ in the sense of 
Definition~\ref{iforc}) for all $i\in\avl\rho$.\qed
\edf 
We shall use Greek characters $\pi,\,\rho,\,\vt$ to denote 
forcing conditions in $\Pi$.

Take notice that if the structures $\gA$ and $\gB$ are not 
elementarily equivalent then $\dP_{\cL_i\,\gA_i\,\gB_i}$ is 
empty; so that in this case $i\not\in\avl \pi$ for all 
$\pi\in\Pi.$ The parameter $w=w_i$ does not participate 
in the definition of $\Pi;$ its role will be to make $\Pi$ 
homogeneous enough to admit a restriction theorem. 


\subsection{The generic extension}
\label{piext}

The aim of this section (Theorem~\ref{conip} below) is to 
prove that \dd\Pi generic extensions of $\dH$ satisfy $\HST+\ip$ 
provided $\dH$ satisfies statement~\ref{gdc} of 
Theorem~\ref{li:hst}. Let us fix a \dd\Pi generic over $\dH$ 
set $G\sq\dH$. 
 
We put $\cN=\cN[\Pi]=\bigcup_{\al\in\Ord}\cN_\al,$ $\cN\sq\dH$ is 
the class of 
\index{class!n@$\cN$ of ``names''}%
all \dd\Pi``names'', defined in $\dH$ as in Subsection~\ref{gen}. 
We introduce, following Subsection~\ref{gen},\its

\ben
\def\theenumi{\arabic{enumi})}
\def\labelenumi{\theenumi}

\item a set $a[G]$ for each ``name'' $a\in\cN$ by induction 
on $\nrk a$,\hfill and 
\index{zzag@$a[G]$}
\its

\item the extension $\dH[G]=\ans{a[G]:a\in\cN}$ with  
the membership $\ing$. 
\index{zzhg@$\dH[G]$}
\index{class!hg@$\dH[G],$ generic extension}
\index{membership relation!ing@$\ing$ on $\dH[G]$}
\een

\bdf
\label{ips}
$\ips$ is the strong form of $\ip$ which asserts that any two 
\index{isomorphism property!ips@strong, $\ips$}
internally presented elementarily equivalent structures of a 
first--order language containing (standard size)--many symbols, 
are isomorphic via a {\it locally internal\/} (see 
Definition~\ref{loc}) isomorphism.\qed
\edf


\bte
\label{conip}
Suppose that in addition to\/ \ref{dhis} and\/ \ref{dhwf} 
$\dH$ satisfies statement~\ref{gdc} in 
Theorem~\ref{li:hst}. Let\/ $\Pi$ be defined as above, in\/ 
$\dH.$ Then every\/ \dd\Pi generic extension\/ $\dH[G]$ is a model 
of\/ $\HST,$ a plain extension of\/ $\dH,$ where the ``strong'' 
isomorphism property\/ $\ips$ holds. 
\ete
This is the main result of this section. 
We begin the proof with several introductory remarks mainly 
devoted to relationships between the model $\dH[G]$ and its 
submodels.

We observe that $\Pi$ is a proper class, not a set in 
$\dH.$ This makes it necessary to change something in the 
reasoning in Section~\ref{f}. 
For instance now $\brG=\ans{\ang{\pi,\brpi}:\pi\in\Pi}$ is not a 
set in $\dH,$ so that one cannot assert that $G\in\dH[G].$ However 
this is not a problem because we are now interested in certain 
small parts of $G,$ rather than $G$ itself, to be elements 
of $\dH[G]$. 

Let $C\in\dH,$ $C\sq\Ind.$ Then $\Pi_C=\ans{\pi\in\Pi:\avl\pi\sq C}$ 
\index{zzpic@$\Pi_C$}
\index{zzgc@$G_C$}
is a \underline{set} in $\dH.$ (Use the $\HST$ Collection.) 
We define $G_C=\Pi_C\cap G,$ for each $G\sq\Pi$. 

Let, for $\pi\in\Pi,$ $\pi\res C$ be the restriction of $\pi$ to 
\index{zzpiresc@$\pi\res C$}
the domain $\avl\pi\cap C;$ $\pi\res C\in\Pi$ and $\in \Pi_C$ 
provided $\avl\pi\cap C$ is internal. 
Furthermore we have $G_C=\ans{\pi\res C:\pi\in G}$ provided 
$G\sq\Pi$ is generic and 
$C$ is internal in $\dH$. 


We define, in $\dH,$ a set $\davl a\sq\Ind$ for each ``name'' 
\index{zzadavl@$\davl a$}
$a\in\cN,$ by induction on $\nrk a$ as follows. If $a\in\cN_0$  
then $\davl a=\emps.$ Otherwise we put  
$\davl a=\bigcup_{\ang{\pi,b}\in a}(\davl b\cup\avl\pi).$  
We let $\cN\res C=\ans{a\in\cN:\davl a\sq C},$ for each set 
\index{zznresc@$\cN\res C$}
$C\sq\Ind.$ Then $\cN\res C$ is precisely the class of all 
\dd{\Pi_C}``names''. 

For instance $\brG_C=\ans{\ang{\pi,\brpi}:\pi\in\Pi_C}$ belongs 
to $\cN\res C$. 

\bpro
\label C
Let\/ $G\sq\Pi$ be\/ \dd\Pi generic over $\dH.$ Suppose that\/ 
$C\sq\Ind$ is an internal set. Then\its
\ben
\itla{c1}\msur 
$G_C=\brG_C[G]\in\dH[G]$ is\/ \dd{\Pi_C}generic over\/ $\dH$.\its 

\itla{c2}\msur
$\dH[G_C]=\ans{a[G_C]:a\in\cN\res C}$ is a transitive subclass 
of\/ $\dH[G]$.\its

\itla{c3} 
If\/ $a\in\cN\res C$ then\/ $a[G]=a[G_C]$.
\een
\epro
\proof An ordinary application of the product forcing 
technique.\qed\eproof

It follows that $\dH[G]$ 
is a plain extension of $\dH,$ by Theorem~\ref{is}. 

\subsection{The product forcing relation}
\label{pifo}

The continuation of the proof of Theorem~\ref{conip} involves  
forcing. 

There is a problem related to forcing: the definition of $\fo$ for 
the atomic formulas $a=b$ and $b\in a$ in Subsection~\ref{fo} becomes 
unsound in the case when the notion of forcing is not a set in the 
ground model $\dH,$ as in the case we consider now. (This is a 
problem in the $\ZFC$ setting of forcing as well, see~\cite{sh}.) 
The solution follows the $\ZFC$ patterns: the inductive definition 
of forcing for atomic formulas can be executed using only 
set parts $\Pi_C$ of the whole forcing $\Pi$. 

For each set $C\sq\Ind,$ let ${\spfo C}$ and $\pfo C$ be the 
\index{forcing!pfo@$\pfo C$}
\index{forcing!spfo@$\spfo C$}
forcing relations associated in $\dH$ with $\Pi_C$ as the forcing 
notion, as in Subsection~\ref{fo}. 

Our plan is as follows. We define the \dd\Pi forcing $\fo$ for 
atomic formulas $a=b$ and $b\in a$ ($a,\,b$ being ``names'' in 
$\cN$) using $\pfo C$ for sufficiently large internal sets 
$C\sq\Ind.$ Then we define $\fo$ for other formulas following 
the general construction (items \ref{f-st} through \ref{f-all} in 
Subsection~\ref{fo}). 
%

To start with, let us describe connections between $\spfo C$ 
for different $C,$ and the relation $\pi\sfo b\in a$ defined 
\index{forcing!sfo@$\sfo$}
for $\Pi$ as the notion of forcing as in Subsection~\ref{fo}. 

\bpro
\label{.s}
Let\/ $a,\,b\in\cN,\;\,\pi\in\Pi.$  Suppose that\/ $C$ is an 
internal set, and\/ $\davl a\cup\davl b\sq C\sq \Ind.$ 
Then\/ $\pi\sfo b\in a$ iff\/ $\pi\res C\spfo C b\in a$.
\epro
\proof Elementary verification, based 
on the fact that $\ang{\pi,b}\in a$ implies 
$\avl\pi\sq\davl a,$ is left for the reader.\qed\eproof

The following lemma is of crucial importance. 

\ble
\label{.}
Let\/ $\Phi$ be a formula of the form\/ $a=b$ or\/ $b\in a,$ 
where\/ $a,\,b\in\cN.$ Suppose that\/ $C,\,C'$ are internal 
sets, and\/ $\davl a\cup\davl b\sq C\sq C'\sq\Ind.$ Let 
finally\/ $\pi'\in\Pi_{C'}$ and\/ $\pi=\pi'\res C.$ Then\/ 
$\pi'\pfo{C'}\Phi$ iff\/ $\pi\pfo C\Phi$.
\ele
\proof The proof goes on by induction on the 
ranks $\nrk a$ and $\nrk b$.

Let $\Phi$ be the formula $b=a.$ Let us suppose that 
$\pi'\pfo{C'}b=a,$ and prove $\pi\pfo C b=a.$ Let 
$\rho\in\Pi_C,$ $\rho\<\pi,$ and $\rho\spfo{C}x\in a.$ We define  
$\rho'=\rho\cup (\pi'\res(C'\setminus C))\in\Pi_{C'};$ then 
$\rho'\res C=\rho,$ and $\rho'\<\pi'$ because 
$\rho\<\pi.$ Take notice that $\davl x\sq\davl a\sq C,$ so we have 
$\rho'\spfo{C'}x\in a$ by Proposition~\ref{.s}. It follows that 
$\rho'\pfo{C'}x\in b,$ since $\pi'\pfo{C'}b=a$ was assumed. We 
finally have $\rho\pfo{C}x\in b$ by the induction hypothesis.

Conversely, suppose that $\pi\pfo C b=a$ and prove 
$\pi'\pfo{C'}b=a.$ Assume that $\rho'\in\Pi_{C'},$ $\rho'\<\pi',$ 
and $\rho'\spfo{C'}x\in a.$ Then $\davl x\sq\davl a\sq C,$ so that 
$\rho=\rho'\res C$ satisfies $\rho\<\pi$ and $\rho\spfo{C}x\in a$ 
by Proposition~\ref{.s}. Since $\pi\pfo C b=a,$ we have 
$\rho\pfo{C} x\in b.$ We conclude that $\rho'\pfo{C'} x\in b$ 
by the induction hypothesis. 

Let $\Phi$ be the formula $b\in a.$ Suppose that 
$\pi'\pfo{C'}b\in a$ and prove $\pi\pfo C b\in a.$ Let 
$\rho\<\pi.$ Then 
$\rho'=\rho\cup (\pi'\res(C'\setminus C))\in\Pi_{C'},$ 
$\rho'\<\pi',$ so that there exist $\vt'\in\Pi_{C'},$ 
$\vt'\<\rho',$ and some $z$ such that $\vt'\spfo{C'}z\in a$ 
and $\vt'\pfo{C'}b=z.$ Then $\vt=\vt'\res C\in\Pi_C$ and 
$\vt\<\rho.$ On the other hand, $\davl z\sq\davl a\sq C,$ so that 
$\vt\spfo{C} z\in a$ and $\vt\pfo C b=z,$ as required. 

Conversely, suppose that $\pi\pfo{C}b\in a$ and prove 
$\pi'\pfo{C'} b\in a.$ Let $\rho'\in\Pi_{C'},$ $\rho'\<\pi'.$ 
Then $\rho=\rho'\res C\in\Pi_C$ and $\rho\<\pi,$ so that  
$\vt\spfo{C} z\in a$ and $\vt\pfo C b=z,$ for some 
$\vt\in\Pi_C,$ $\vt\<\rho,$ and a ``name'' $z.$ We put 
$\vt'=\vt\cup (\rho'\res(C'\setminus C));$ so that 
$\vt'\in\Pi_{C'},$ $\vt'\<\rho',$ and $\vt=\vt'\res C.$ 
Then $\vt'\spfo{C'} z\in a,$ and $\vt'\pfo{C'} b=z$ by the 
induction hypothesis.\qed\eproof

This leads to the following definition. 

\bdf
\label{Pi}
We introduce the forcing relation ${\fo}={\fo_\Pi}$ as follows.
\index{forcing!fo@$\fo$} 

Let $a,\,b\in\cN$ and $\pi\in\Pi.$ We set $\pi\fo b\in a$ iff 
$\pi\pfo{C}b\in a,$ and $\pi\fo b=a$ iff $\pi\pfo{C}b=a,$ 
whenever $C\sq\Ind$ is an internal set satisfying 
$\avl\pi\cup\davl a\cup\davl b\sq C.$ 
(This does not depend on the choice of $C$ by Lemma~\ref..)

The relation $\fo$ expands on the standardness predicate and 
non--atomic formulas in accordance with items \ref{f-st} through 
\ref{f-all} of Subsection~\ref{fo}.\qed
\edf
In view of this definition, Lemma~\ref. takes the following 
form:

\bcor
\label{;}
Let\/ $\Phi$ be a formula of the form\/ $a=b$ or\/ 
$b\in a,$ where\/ $a,\,b$ are names in $\cN.$ Suppose 
that\/ $\pi\in\Pi,$ and\/ $C\sq\Ind$ is an internal set, 
satisfying\/ $\davl a\cup\davl b\sq C.$ Then\/ $\pi\fo\Phi$ 
iff\/ $\pi\res C\pfo C\Phi$.\qed
\ecor
Furthermore, it occurs that the forcing $\fo$ still obeys the 
general scheme !

\bcor
\label{same}
The relation\/ $\fo$ formally satisfies requirements of items\/ 
\hbox{\ref{f=}} and\/ \hbox{\ref{f-in}} of the definition of 
forcing in Subsection~\ref{fo}, with respect to\/ $\Pi$ as the 
forcing notion. 
\ecor
\proof An easy verification, with reference to 
Corollary~\ref; for large enough internal sets $C$.\qed

\brem
\label{rems}
Corollary~\ref{same} guarantees that the results 
obtained for ``set'' size forcing in subsections \ref{fo} and 
\ref{trul} remain valid for the forcing $\fo$ associated with 
$\Pi,$ with more or less the same proofs.\qed
\erem 

We need to verify two particular properties of the notion of 
forcing $\Pi,$ before the proof of Theorem~\ref{conip} starts. 
One of them deals with the restriction property of the forcing: 
we would like to prove that $p\fo \Phi$ iff $p\res C\fo\Phi$ 
provided $\davl\Phi\sq C,$ for all, not only atomic formulas 
$\Phi.$ The other one is the standard size distributivity of 
$\Pi$.  

\subsection{Automorphisms and the restriction property}
\label{aut}

We apply a system of automorphisms of the notion of forcing $\Pi$ 
to approach the restriction property.

Let $D\sq\Ind$ be an internal set. An internal bijection 
$h:D\,$onto$\,D$ satisfying the requirement:\its
\ben
\def\theenumi{$(\fnsymbol{enumi})$} 
\def\labelenumi{\theenumi}
\itla{yy} 
If $i=\ang{w,\kpa,\fL,\fA,\fB}\in D$ then 
$h(i)=\ang{w',\kpa,\fL,\fA,\fB}$ for some (internal) $w'$ 
and the same $\kpa,\,\fL,\,\fA,\,\fB$,\its
%
\een
will be called a {\it correct\/} bijection. 
\index{correct bijection}
In this case we define $H(i)=h(i)$ for $i\in D,$ and 
$H(i)=i$ for $i\in\Ind\setminus D;$ so that $H=H_h$ $1-1$ maps 
\index{zzhh@$H_h$}
$\Ind$ onto $\Ind.$ $H$ obviously inherits property \ref{yy}. 

The bijection $H$ generates an order automorphism of the notion 
of forcing $\Pi, $ defined as follows. Let $\pi\in\Pi.$ We define 
$H\pi\in\Pi$ so that $\avl{H\pi}=\ans{H(i):i\in\avl\pi}$ and 
\index{zzhpi@$H\pi$}
$(H\pi)_{H(i)}=\pi_i$ for each $i\in\avl\pi.$ It follows from 
\ref{yy} that the map ${\pi\map H\pi}$ is an order 
automorphism of $\Pi$. 

Let us expand the action of $H$ onto ``names''. We define, in 
$\dH,$ $H[a]$ for each ``name'' $a,$ by induction on $\nrk a.$ If 
\index{zzha@$H[a]$} 
$a=\brx\in\cN_0$ then we put $H[a]=a.$ If $\nrk a>0$ then we 
set $H[a]=\ans{\ang{H\pi,H[b]}:\ang{\pi,b}\in a}.$ One easily proves 
that $H[a]\in\cN$ and $\nrk a=\nrk H[a]$.

For a \ste-formula $\Phi$ containing ``names'' in $\cN,$  
we let $H\Phi$ denote the formula obtained by changing 
each ``name'' $a$ in $\Phi$ to $H[a]$.

\bpro
\label{inv}
Let\/ $h$ be a correct bijection, and\/ $H=H_h.$ For any 
condition\/ $\pi\in\Pi$ and any closed formula\/ $\Phi,$ 
having ``names'' in\/ $\cN$ as parameters, $\pi\fo\Phi$ iff\/ 
$H\pi\fo H\Phi$.
\epro
\proof  We omit the routine verification, which can be conducted 
by induction on the complexity of the formulas involved, 
following arguments known from the theory of generic extensions 
of models of $\ZFC$.\qed

\bcor
\label{rtion}
{\rm (Restriction)} \ 
Suppose that\/ $\pi\in\Pi,$ $\Phi$ is a closed formula containing 
``names'' in\/ $\cN$ as parameters, and\/ $\pi\fo\Phi.$ Suppose 
also that\/ $C$ is an internal set, and $\davl\Phi\sq C.$ Then 
$\pi\res C\fo\Phi$.
\ecor
\proof It follows from Lemma~\ref{ok} that otherwise there 
exists a pair of conditions $\pi,\,\rho\in\Pi$ such that 
$\pi\res C=\rho\res C,$ $\pi\fo\Phi,$ but $\rho\fo\neg\;\Phi.$ 
Let $D=\avl\pi,$ $E=\avl{\rho}.$ It is clear that there exists 
an internal set $W$ satisfying $C\cup D\cup E\sq W,$ and an 
internal {\it correct\/} bijection $h:W\,\hbox{ onto }\,W$ which 
is the identity on $C$ and satisfies $E\cap (h\ima D)\sq C.$ Let 
$H=H_h$ be defined, from $h,$ as above. Let $\pi'=H\pi.$ Then 
$\pi'\res C=\pi\res C=\rho\res C$ because $h\res C$ is the 
identity. Furthermore $\avl{\pi'}=h\ima D,$ so that 
$\avl{\pi'}\cap\avl\rho\sq C.$ We conclude that $\pi'$ and 
$\rho$ are compatible in $\Pi$. 

On the other hand, $\pi'\fo H\Phi$ by Proposition~\ref{inv}. Thus 
it suffices to demonstrate that $\Phi$ coincides with $H\Phi.$ 
We recall that $\davl\Phi\sq C,$ so that each ``name'' $a$ which 
occurs in $\Phi$ satisfies $\davl a\sq C.$ However one can easily 
prove, by induction on $\nrk a$ in $\dH,$ that $H[a]=a$ 
whenever $\davl a\sq C,$ using 
the fact that $h\res C$ is the identity. We conclude that 
$H\Phi$ is $\Phi,$ as required.\qed

\subsection{Standard size distrubutivity of the product forcing}
\label{picompl}

We are going to prove that $\Pi$ is standard size distributive 
in $\dH$ provided $\dH$ satisfies requirement~\ref{gdc} of 
Theorem~\ref{li:hst}. 

\bpro
\label{tfo1}
$\Pi$ is standard size closed in $\dH$.
\epro 
\proof In principle the proof copies that of Proposition~\ref{fo1}, 
but we need to take more time to reduce the problem to Saturation. 
Suppose that $\la$ is a cardinal, $\pi_\al\;(\al<\la)$ are 
conditions in $\Pi,$ and $\pi_\ba\<\pi_\al$ 
whenever $\al<\ba<\la.$ Using the $\HST$ Collection and 
Lemma~\ref{hstb} in $\dH,$ we get a standard set $S\sq\Ind$ such 
that each $\pi_\al$ in fact belongs to $\Pi_S.$ Let us check that 
the sequence has lower bound already in $\Pi_S$.

We observe 
that by the Collection axiom again, there exists a cardinal $\kpa$ 
such that $\kpa_i\<\upa\kpa$ in $\dI$ whenever $i\in S.$ ($\kpa_i$ 
was defined in Subsection~\ref{prod}.) Then $\Pi_S$ is an 
intersection of \dd{(\<\!\kpa)}-many internal sets by definition. 
Therefore every set  
$P_\al=\ans{\pi\in\Pi_S:\pi\<\pi_\al}\;\;(\al<\la)$ is, uniformly 
on $\al<\la,$ an intersection of \dd{(\<\!\kpa)}-many internal 
sets, too. Furthermore the sets $P_\al$ are nonempty and 
$P_\ba\sq P_\al$ whenever $\al<\ba<\la.$ Since $\la$ and $\kpa$ are 
sets of standard size by Lemma~\ref{wo=ss}, we can use Saturation 
to obtain $\bigcap_{\al<\la}P_\al\not=\emps,$ as required.\qed

\bpro
\label{tgdc}
Assume that the ground model\/ $\dH$ satisfies 
statement~\ref{gdc} of Theorem~\ref{li:hst}. Then\/ 
$\Pi$ is standard size distributive in\/ $\dH$.
\epro
\proof We cannot directly refer to \ref{gdc} and 
Proposition~\ref{tfo1} because $\Pi$ is a proper class rather than 
a set in $\dH.$ (Being a set is essential in the proof of 
Theorem~\ref{li:hst}, by the way.) But of course we shall 
reduce the problem to the assumption of \ref{gdc}, by the choice 
of a suitable set part of $\Pi$. 

We have to prove the following. Let $\kpa$ be a cardinal in $\dH$ 
(in the sense of Subsection~\ref{orcar}), and $D$ be a \ste-definable 
in $\dH$ subclass of $\kpa\ti\Pi.$ Suppose that each class 
$D_\al=\ans{\pi:\ang{\al,\pi}\in D}$ is open dense in $\Pi.$ Then 
the intersection $\bigcap_{\al<\kpa}D_\al$ is dense in $\Pi$ as 
well.

To prove the assertion, we fix a condition $\pi_0\in\Pi.$ Let 
$S_0\sq\Ind$ be an arbitrary standard set such that 
$\pi_0\in \Pi_{S_0}.$ Let $\kpa^+$ be the next cardinal. 
(In concern of cardinals, we are in $\dV,$ a $\ZFC$ universe.) 
Let us define an increasing sequence of standard sets 
$S_\al\sq\Ind\;\,(\al<\kpa^+)$ as follows. 
$S_0$ already exists. Suppose that $\ga<\kpa$ and $S_\al$ is 
defined for each $\al<\ga.$ We first put 
$S'_\ga=\bigcup_{\al<\ga} S_\al.$ (For instance $S'_\ga=S_\ba$ 
provided $\ga=\ba+1.$) Using Collection and Lemma~\ref{hstb} in 
$\dH,$ and the assumption that every $D_\al$ is dense, we obtain 
a standard set $S,\,\; S'_\ga\sq S\sq\Ind,$ such that for any 
$\pi\in\Pi_{S'_\ga}$ ($\Pi_{S'_\ga}$ is a set !) and any $\al<\kpa$ 
there exists 
$\rho\in\Pi_S\cap D_\al$ such that $\rho\<\pi$ in $\Pi.$ Let 
$S_\ga$ denote the least standard set $S$ of the form 
$V_\nu\cap\Ind$ (where $V_\nu$ is the \dd\nu th level of the 
von Neumann hierarchy; $\nu$ being an \dd\dS ordinal) satisfying 
this property. 

Let $S=\bigcup_{\al<\kpa^+} S_\al.$ Then $P=\Pi_S$ is a 
{\it set\/}. Furthermore $\pi_0\in P,$ and each intersection 
$D'_\al=D_\al\cap P$ is dense in $P$ by the construction, 
because $P=\bigcup_{\al<\kpa^+}\Pi_{S_\al}.$ (In this argument, 
we use Saturation and the fact that $\avl\pi$ is internal for 
$\pi\in\Pi.$) It remains to check that $P$ 
is \dd\kpa closed in $\dH:$ every decreasing sequence 
$\ang{\pi_\al:\al<\kpa}$ has a lower bound in $P.$ 
(Indeed then $P$ is \dd\kpa distriburive by the assumption of 
\ref{gdc}, so that the intersection $\bigcap_{\al<\kpa}D'_\al$ 
is dense in $P,$ etc.) 

By Proposition~\ref{tfo1}, the sequence has a lower bound 
$\pi\in \Pi.$ (We cannot run the proof of Proposition~\ref{tfo1} 
for $P$ directly because $P$ is not a standard size intersection 
of internal sets.) Since the construction of 
$S_\al$ involves all ordinals $\al<\kpa^+,$ there exists an ordinal 
$\ga<\kpa^+$ such that every condition $\pi_\al\;\,(\al<\kpa)$ 
belongs to $\Pi_{S_\ga}.$ Then $\rho=\pi\res S_{\ga+1}$ still 
satisfies $\rho\<\pi_\al$ for all $\al,$ but $\rho\in P,$ 
as required.\qed

\subsection{Verification of the axioms}

This subsection starts the proof of Theorem~\ref{conip}. 

The verification of $\HST$ in $\dH[G]$ copies, to some extent, 
the proof of Theorem~\ref{t:hg}. (The standard size distributivity  
of the forcing, assumed in Theorem~\ref{t:hg}, now follows from 
Proposition~\ref{tgdc}.) Only the proofs of Separation and 
Collection need to be performed anew, because it was essential  
in Subsection~\ref{hg:hst} that the notion of 
forcing is a set in the ground model. 

{\it Separation\/}. We follow the reasoning in the proof of 
Theorem~\ref{t:hg}. Suppose that $X\in\cN,$ and $\Phi(x)$ is a 
\ste-formula which may contain ``names'' in $\cN$ as parameters. 
We have to find a ``name'' $Y\in \cN$ satisfying the equality 
$Y[G]=\ans{x\in X[G]:\Phi[G](x)}$ in $\dH[G]$. 

We observe that all elements of $X[G]$ in $\dH[G]$ 
are of the form $x[G]$ where $x$ belongs to the set of ``names'' 
$\cX=\ans{x\in\cN:\exists\,\pi\;(\ang{\pi,x}\in X)}\in\dH.$  
We cannot now define 
$Y=\ans{\ang{\pi,x}\in \Pi\ti \cX:\pi\fo x\in X\cj\Phi(x)},$ 
simply because this may be not a set in $\dH.$ (We recall that 
$\Pi$ is a proper class in $\dH$.) To overcome this difficulty, 
we replace $\Pi,$ using the restriction theorem, by a suitable 
$\Pi_C$.

It follows from Lemma~\ref{hstb} that there exists an internal 
(even standard) set $C\sq\Ind$ such that 
$\davl\Phi\cup\davl X\sq C.$ Then, by the way, $\davl x\sq C$ 
for every $x\in \cX.$ We set 
${Y=\ans{\ang{\pi,x}\in\Pi_C\ti \cX:\pi\fo x\in X\cj\Phi(x)}}.$ 
One easily proves that $Y$ is the required ``name'', using 
Corollary~\ref{rtion} (the restriction theorem) and following 
usual patterns. (See e.\ g.\ Shoenfield~\cite{sh}.) 

{\it Collection\/}. We suppose that $X\in\cN,$ and $\Phi(x,y)$ is a 
formula  with ``names'' in $\cN$ as parameters. Let $\cX\sq\cN,$ 
$\cX\in\dH,$ be 
defined in $\dH$ as in the proof of Separation. It would suffice to 
find a set of ``names'' $\cY\in\dH,$ $\cY\sq\cN,$ such that for 
every $x\in\cX$ and every condition $\vt\in\Pi,$ if 
$\vt\fo\exists\,y\;\Phi(x,y)$ then there exist: a ``name'' $y\in\cY$ 
and a stronger condition $\rho\<\vt$ which forces $\Phi(x,y)$. 

Let us choose an internal set $C_0$ so that 
$\davl\Phi\cup\davl X\sq C_0$.

We have to be careful because $\Pi,$ the notion of forcing, is a 
proper class in $\dH.$ However, since $\Pi_{C_0}$ is a set in $\dH,$ 
there exist: a set $P\sq\Pi$ of forcing conditions, and a set 
$\cY_0\in\dH,$ $\cY_0\sq\cN,$ of ``names'', satisfying the property: 
if $x\in \cX,$ and $\pi_0\in\Pi_{C_0}$ forces $\exists\,y\;\Phi(x,y)$ 
then there exist: a condition $\pi\in P,$ $\pi\<\pi_0,$ and a 
``name'' $y\in\cY_0,$ such that $\pi\fo\Phi(x,y)$. 

The set $\cY_0$ is not yet the $\cY$ we are looking for. To get 
$\cY,$ we first of all choose an internal set $C$ such that 
$C_0\sq C,$ $\avl\pi\sq C$ for all $\pi\in P,$ $\davl y\sq C$ for all 
$y\in \cY_0,$ the difference $C\setminus C_0$ is \dd\dI infinite, and 
moreover, for any $i=\ang{w,\kpa,\fL,\fA,\fB}\in C$ there exist 
\dd\dI infinitely many different indices $i'\in C$ of the form 
$i'=\ang{w',\kpa,\fL,\fA,\fB}\in C$ (with $w\not=w'$ but the same 
$\kpa,\,\fL,\,\fA,\,\fB$). 
Each internal {\it correct\/} bijection $h:C\,\hbox{ onto }\,C$ 
generates an automorphism $H_h$ of $\Pi,$ see Subsection~\ref{aut}. 
Let us prove that 
\dm
\cY=\ans{H_h[y]:y\in \cY\;\hbox{ and }\;h\in\dI\;
\hbox{ is a {\it correct\/} bijection }\,C\,\hbox{ onto }\,C}
\dm
is a set of ``names'' satisfying the property we need. 
(To see that $\cY$ is a set in 
$\dH$ notice the following: all the bijections $h$ considered are 
internal by definition, so we can use internal power sets in $\dI$.)

Let $x\in \cX$ and $\vt\in\Pi.$ 
Suppose that $\vt\fo\exists\,y\;\Phi(x,y).$ Then the condition 
$\pi_0=\vt\res C_0$ also forces $\exists\,y\;\Phi(x,y)$ by 
Corollary~\ref{rtion}. ($\davl{\exists\,y\;\Phi(x,y)}\sq C_0$ by 
the choice of $x$ and $C_0$.) Then, by the choice of $P$ and $\cY_0,$ 
there exist: a condition $\pi\in P,$ $\pi\<\pi_0,$ and a ``name'' 
$y\in \cY_0,$ such that $\pi\fo\Phi(x,y).$ Unfortunately $\pi$ 
may be incompatible with $\vt;$ otherwise we would immediately 
consider any condition $\rho$ stronger than both $\pi$ and $\vt.$ 
To overcome this obstacle, let us use 
an argument from the proof of Corollary~\ref{rtion}. 

Let $\vt'=\vt\res C.$ Take notice that $E=\avl{\pi}$ and 
$D'=\avl{\vt'}$ are, by definition, \dd\dI{\it finite\/}~\footnote 
{\ This is the only point where the finiteness of the domains 
$\avl\pi,$ $\pi\in \Pi,$ see Definition~\ref{totp}, is used. In 
fact the proof does not change much if the domains $\avl\pi$ are 
restricted to be less than a fixed \dd\dI cardinal.} internal subsets of $C.$   
There exists, by the choice of $C,$ an internal {\it correct\/} 
bijection $h:C\,\hbox{ onto }\,C$ such that $h\res C_0$ is the 
identity and $(h\ima E)\cap D'\sq C_0.$ Let $H=H_h.$ Then 
$\pi'=H\pi\in\Pi_{C},$ $\pi'\res C_0=\pi\res C_0\<\pi_0,$ 
and $\avl{\pi'}\cap\avl{\vt'}\sq C_0,$ so that $\vt'$ and $\pi'$ 
are compatible. Therefore $\pi'$ is also compatible with $\vt$ 
because $\pi'\in\Pi_{C}$ and $\vt'=\vt\res C.$ Let $\rho\in\Pi$ be 
a condition stronger than both $\pi'$ and $\vt$. 

We observe that $\pi'\fo H\Phi(H[x],H[y]),$ by Theorem~\ref{inv}. 
But, $\davl\Phi\sq C_0$ and $\davl x\sq C_0$ by the choice of $C_0,$ 
so that $H\Phi$ coincides with $\Phi$ and $H[x]=x$ because  
$H\res C_0$ is the identity. We conclude that $\rho\fo\Phi(x,y'),$ 
where $y'=H[y]$ is a ``name'' in $\cY$ by definition, as required.

\subsection{Verification of the isomorphism property in the 
extension}

We accomplish the proof of Theorem~\ref{conip} in this subsection.

Since the standard sets are essentially the same in $\dH$ and 
$\dH[G],$ the condensed subuniverse $\dV$ is also one and the 
same in the two universes. Therefore $\dH[G]$ contains the same 
ordinals as $\dH$ does. (See subsections \ref{conden} and 
\ref{orcar}.) Since standard size subsets of $\dH$ in $\dH[G]$ 
all belong to $\dH,$ cardinals in $\dH[G]$ are the same as 
cardinals in $\dH.$ (We recall that by definition 
{\it cardinals\/} mean: well-orderable cardinals, in $\HST$.) 

This reasoning shows that all the triples: 
{\it language -- structure -- structure\/}, to be considered in 
the frameworks of the isomorphism property in $\dH[G],$ are 
already in $\dH.$ Thus let $\cL\in\dH$ be a standard 
size first--order language, containing $\kpa$ symbols in 
$\dH$ ($\kpa$ being a cardinal in $\dH$), and $\gA,\,\gB$ be a 
pair of internally presented \dd{\cL}structures in $\dH.$ We 
have to prove that $\gA$ is isomorphic to $\gB$ in $\dH[G]$.

Using Lemma~\ref{exten} in $\dH,$ we obtain an internal 
first--order language $\fL=\ans{s_\al:\al<\upa\kpa},$ containing 
$\upa\kpa$ symbols in $\dI,$ and internal \dd\fL structures $\fA$ 
and $\fB,$ such that 
$\cL=\upsg\fL=\ans{s_\al:\al<\upa\kpa\cj\st\al},$ and 
$\gA,\,\gB$ are the corresponding restrictions of $\fA,\,\fB.$ In 
other words, $i=\ang{0,\upa\kpa,\fL,\fA,\fB}$ belongs to $\Ind$ 
and $\cL=\cL_i,$ $\gA=\gA_i,$ $\gB=\gB_i$. 

We observe that the set $G_i=\ans{\pi_i:\pi\in G\cj i\in\avl\pi}$ 
belongs to $\dH[G].$ 
(Indeed, since $\Pi_i=\dP_{\cL_i\,\gA_i\,\gB_i}$ is a set in 
$\dH,$ a ``name'' for $G_i$ can be defined in $\dH$ as the set of 
all pairs $\ang{\pi,p},$ where $p\in\Pi_i$ and 
$\pi=\ang{i,p}\in\Pi$ -- so that $\avl\pi=\ans{i}$ and $\pi_i=p$.) 
Furthermore $G_i$ is \dd{\dP_{\cL_i\,\gA_i\,\gB_i}}generic over 
$\dH.$ (An ordinary product forcing argument.) It follows from 
Theorem~\ref{t:isom} in Section~\ref{isom} that $\gA$ and $\gB$ 
are isomorphic in $\dH[G_i]$ via the locally internal isomorphism 
$F_i=\bigcup G_i,$ therefore in $\dH[G],$ as required.

This ends the proof of Theorem~\ref{conip}.\qed

\np

%
\section{Proof of the main theorem}
\label{final}

In this section, we gather the material of the model constructions 
above, with some results from \cite{hyp1,hyp2,hyp3}, to accomplish 
the proof of the main theorem (Theorem~\ref{maint}). 

Let us fix a countable model $\dS$ of $\ZFC.$ We shall {\it not\/} 
assume that $\dS$ is a transitive model; in particular the 
membership relation $\ins$ acting in $\dS$ may be not equal to 
the restriction ${\in}\res\dS$. 

\bpro
\label{h1}
There exists a countable model\/ $\dI$ of\/ $\BST,$ bounded set 
theory, such that the class of all standard sets in\/ $\dI$ 
coincides with\/ $\dS,$ in particular\/ ${\ins}={{\ini}\res\dS}$.
\epro
\proof We refer to \cite{hyp1}, Theorem~2.4. The proof goes on 
as follows. We first add to $\dS$ a generic global choice 
function $\fG,$  
using the method of Felgner~\cite{fe}. This converts $\dS$ into a 
model $\ang{\dS;\fG}$ of $\ZFC$ plus Global Choice, but with the 
same sets as $\dS$ 
originally had. (The assumption of countability of $\dS$ is used 
to prove the existence of a Felgner--generic extension of~$\dS$.) 

The global choice function makes it possible to define, in 
$\ang{\dS;\fG},$ a certain increasing sequence of 
class--many ``adequate'' ultrafilters. The corresponding 
ultralimit of $\dS$ can be taken as $\dI$.\qed

\bpro
\label{h2}
There exists a countable model\/ $\dH'$ of\/ $\HST,$ such that 
the classes of all standard and internal sets in\/ $\dH'$ 
coincide with resp\/{.}\ $\dS$ and $\dI,$ in particular\/ 
${\ini}={{\inhp}\res\dI}$.
\epro
\proof We refer to \cite{hyp2}, Theorem~4.11. To get $\dH',$ we 
first consider $\dE,$ the class of all elementary external sets 
(i.\ e. \ste-definable in $\dI$ subclasses of sets in $\dI,$ see 
Subsection~\ref{ees} above), then define $\dH$ as the collection of 
all sets obtainable by the assembling construction, described in 
Subsection~\ref{assem}, from wf pairs in $\dE.$ Thus in principle 
the $\dH$ obtained this way is equal to $\dL[\dI],$ but we rather 
put this as a separate step.\qed\eproof

We let $\dH$ be $\dL[\dI],$ formally defined in $\dH'$.

\bcor
\label{-ip}
$\dH$ is a countable model\/ of\/ $\HST,$ such that the classes 
of all standard and internal sets in\/ $\dH$ coincide with 
resp\/{.}\ $\dS$ and $\dI,$ the isomorphism property\/ $\ip$ 
fails, and every standard size closed p.\ o.\ set is standard 
size distributive.
\ecor
\proof We refer to Theorem~\ref{li:hst} above.\qed\eproof

This ends the proof of items \ref{esis}, \ref{serv}, \ref{me} of 
Theorem~\ref{maint}, with respect to the theory $\HST+\neg\;\ip.$ 
Let us consider the other one, $\HST+\ip$.

\bcor
\label{+ip}
There exists a countable model\/ $\dH^+$ of\/ $\HST,$ such that 
the classes of all standard and internal sets in\/ $\dH^+$ coincide 
with resp\/{.}\ $\dS$ and $\dI,$ and the strong isomorphism 
property\/ $\ips$ holds.\qed
\ecor
\proof Let us assume, for a moment, that $\dS,$ the initial model 
of $\ZFC,$ is a wellfounded model, in the sense that the membership 
relation $\ins$ is wellfounded in the wider universe. In this case, 
$\dH$ is wellfounded over $\dI$ in the sense of Section~\ref{f} 
because the ordinals in $\dH$ are the same as in $\dV,$ the 
condensed subuniverse, and therefore order isomorphic to the 
ordinals in $\dS.$ Since $\dH$ is countable, there exists a 
\dd\Pi generic extension $\dH^+=\dH[G],$ of the type considered 
in Section~\ref{total}. $\dH^+$ is the required model by 
Theorem~\ref{conip}. 

Let us now consider the general case: $\dS$ at the beginning, 
and $\dH$ at the end, may be not wellfounded. Then of course 
one cannot carry out the construction of $\dH[G]$ described in 
Subsection~\ref{gen}. 

But one can conduct a different construction, also known from 
manuals on forcing for models of $\ZFC.$ This construction goes 
on as follows. We first define the forcing relation $\fo,$ 
associated with $\Pi$ in $\dH,$ as in Subsection~\ref{pifo}, 
which does not need any previous construction of the extension. 
Then we define, given a generic set $G\sq\Pi,$ the relations: 
$a\eqg b$ iff ${\exists\,\pi\in G\;(\pi\fo a\eqg b)},$ and 
similarly $a\ing b$ and $\stg a,$ for all ``names'' 
$a,\,b\in\cN.$ $\eqg$ can be easily proved to be an equivalence 
relation on $\cN,$ while the other two relations to be 
\dd{\eqg}invariant. This allows to define $\dH[G]$ to be the 
quotient $\cN/{\eqg},$ equipped with the quotients of $\ing$ and 
$\stg$ as the atomic relations. The map $x\map\,$(the 
\dd{\eqg}class of $x\,)$ is a \ste-isomorphism $\dH$ onto an 
\dd{\ing}transitive part of $\dH[G].$ (We refer to 
Shoenfield~\cite{sh}.) 

This approach makes it possible to carry out the whole system 
of reasoning used to prove Theorem~\ref{conip}, with minor 
changes.\qed\eproof 

Corollary~\ref{+ip} implies the statements of items \ref{esis}, 
\ref{serv}, \ref{me} of Theorem~\ref{maint}, with respect to the 
theory $\HST+\ip$. 

Let us finally demonstrate that the defined above models,  
$\dH$ of the theory $\HST+\neg\;\ip$ and $\dH^+$ of the theory 
$\HST+\ip,$ satisfy the additional requirement 
\ref{red} of Theorem~\ref{maint}. 

We fix a \ste-formula $\Phi(x_1,...,x_n)$. 

{\it Step 1\/}. Let $\Phi_1(x_1,...,x_n)$ be the formula 
$\emps\fo\Phi(\brx_1,...,\brx_n),$ where $\fo$ is the forcing 
relation $\pfo\Pi,$ associated with $\Pi$ in $\dH,$ while 
$\emps$ is the empty set considered as a forcing condition. 
It is an easy consequence of the restriction theorem 
(Corollary~\ref{rtion}) and the truth lemma 
(Theorem~\ref{truth}) that
\dm
\dH^+\mo\Phi(x_1,...,x_n)
\hspace{5mm}\hbox{iff}\hspace{5mm}
\dH\mo\Phi_1(x_1,...,x_n)\eqno{(1)}
\dm
--- for all $x_1,...,x_n\in\dH.$ This reasoning obviously 
eliminates $\dH^+$ from the problem of consideration, and reduces 
the question to $\dH$. 

{\it Step 2\/}. We recall that, by the construction, $\dH$ is 
$\dL[\dI]$ in a model $\dH'$ of $\HST.$ (In fact $\dH=\dH',$ but 
we shall not use this.) It follows from Proposition~\ref{respoc}, 
applied in $\dH',$ that $\dH$ has a definable interpretation in 
$\dE,$ the collection of all elementary external sets. Therefore 
for each \ste-formula $\Phi_1(x_1,...,x_n)$ there exists another 
\ste-formula $\Phi_2(x_1,...,x_n)$ such that, for all 
$x_1,...,x_n\in\dE$,
\dm
\dH\mo\Phi_1(x_1,...,x_n)
\hspace{5mm}\hbox{iff}\hspace{5mm}
\dE\mo\Phi_2(x_1,...,x_n)\,.\eqno{(2)}
\dm

{\it Step 3\/}. By definition sets in $\dE$ admit a uniform 
\ste-definition from the point of view of $\dI.$ This makes it 
possible to pull things down to $\dI:$ 
for each \ste-formula $\Phi_2(x_1,...,x_n)$ there exists another 
\ste-formula $\Phi_3(x_1,...,x_n)$ such that, 
for all $x_1,...,x_n\in\dI$,
\dm
\dE\mo\Phi_2(x_1,...,x_n)
\hspace{5mm}\hbox{iff}\hspace{5mm}
\dI\mo\Phi_3(x_1,...,x_n)\,.\eqno{(3)}
\dm
 
{\it Step 4\/}. We finally observe that $\dI$ admits a reduction 
to $\dS,$ by a result proved in \cite{hyp1} (Corollary 1.6 there), 
so that for each \ste-formula $\Phi_3(x_1,...,x_n)$ there exists 
a \mem-formula $\Phi_4(x_1,...,x_n)$ such that
\dm
\dI\mo\Phi_3(x_1,...,x_n)
\hspace{5mm}\hbox{iff}\hspace{5mm}
\dS\mo\Phi_4(x_1,...,x_n)\eqno{(4)}
\dm
holds for all $x_1,...,x_n\in\dS$.

Taking the statements $(1)$ through $(4)$ together we conclude 
that the models $\dH$ and $\dH^+$ satisfy the additional requirement 
\ref{red} of Theorem~\ref{maint}.\qed

\np

%

%
\addcontentsline{toc}{section}{\protect\numberline{}Index}
\small
\printindex

\end{document}